\documentclass{article}
\usepackage{amsmath,amssymb}

\usepackage{graphicx,color}
\usepackage{psfrag}
\usepackage{listings}
\usepackage[colorlinks=true,linkcolor=black,citecolor=black,filecolor=blue,urlcolor=blue]{hyperref}

\usepackage[outerbars]{changebar}

\newcommand{\tF}{\tilde F}
\newcommand\DtDx{\frac{\Delta t}{\Delta x}}

\newenvironment{mat}{\left[ \begin{array}{ccccccccccccccc}}{\end{array}\right]}
\newenvironment{rmat}{\left[ \begin{array}{rrrrrrrrrrrrr}}{\end{array}\right]}
\def\bcm{\begin{mat}}   
\def\ecm{\end{mat}}
\def\brm{\begin{rmat}}  
\def\erm{\end{rmat}}

\newenvironment{pvect}{\left( \begin{array}{c}}{\end{array}\right)}
\def\bpvect{\begin{pvect}}
\def\epvect{\end{pvect}}

\newenvironment{pwdef}{\left\{ \begin{array}{ll}}{\end{array}\right.}
\newcommand\bpwdef{\begin{pwdef}}
\newcommand\epwdef{\end{pwdef}}

\def\bsplit{\begin{split}}
\def\esplit{\end{split}}

\newcommand{\Sec}[1]{Section~\ref{sec:#1}}

\newcommand{\ignore}[1]{}
\newcommand{\Comment}[1]{}
\newcommand{\eqn}[1]{(\ref{#1})}
\newcommand{\eq}{\begin{equation}}
\newcommand{\en}{\end{equation}}
\newcommand{\eqm}{\begin{eqnarray}}
\newcommand{\enm}{\end{eqnarray}}
\newcommand{\eqmno}{\begin{eqnarray*}}
\newcommand{\enmno}{\end{eqnarray*}}
\newcommand{\eqml}[1]{\eql{#1}\begin{array}{rcl}}
\newcommand{\enml}{\end{array}\en}
\newcommand{\eql}{\begin{equation}\label}
\newcommand{\eqsub}[1]{\begin{subequations}\label{#1}\eqm }
\newcommand{\ensub}{\enm\end{subequations}}
\newcommand{\half}{{\frac{1}{2}}}

\newcommand{\goto}{\rightarrow}

\def\bc{\begin{center}}
\def\ec{\end{center}}
\def\bi{\begin{itemize}}
\def\ei{\end{itemize}}
\def\be{\begin{enumerate}}
\def\ee{\end{enumerate}}

\def\lam{\lambda}
\def\reals{{{\rm l} \kern -.15em {\rm R} }}

\def\qquad{\quad\quad}

\def\u{\upsilon}

\def\ddt{{\partial\over\partial t}}
\def\ddx{{\partial\over\partial x}}
\def\ddy{{\partial\over\partial y}}
\def\ddz{{\partial\over\partial z}}

\def\inv{^{-1}}

\def\Dx{\Delta x}

\def\Dt{\Delta t}

\newcommand{\imh}{_{i-1/2}}

\newcommand{\iph}{_{i+1/2}}

\newcommand\W{{\cal W}}
\newcommand\tW{\widetilde {\cal W}}

\newcommand\Dq{\Delta Q}

\ignore{

\newcommand{\ignore}[1]{}

\newenvironment{mat}{\left[ \begin{array}{ccccccccccccc}}{\end{array}\right]}
\newenvironment{rmat}{\left[ \begin{array}{rrrrrrrrrrrrr}}{\end{array}\right]}
\newenvironment{lmat}{\left[ \begin{array}{lllllllllllll}}{\end{array}\right]}
\newcommand\bcm{\begin{mat}}
\newcommand\ecm{\end{mat}}
\newcommand\brm{\begin{rmat}}
\newcommand\erm{\end{rmat}}
\newcommand\blm{\begin{lmat}}
\newcommand\elm{\end{lmat}}

\newcommand\bc{\begin{center}}
\newcommand\ec{\end{center}}
\newcommand\bi{\begin{itemize}}
\newcommand\ei{\end{itemize}}
\newcommand\be{\begin{enumerate}}
\newcommand\ee{\end{enumerate}}

\newcommand\bsplit{\begin{split}}
\newcommand\esplit{\end{split}}
\newcommand{\eqn}[1]{(\ref{#1})}
\newcommand{\eq}{\begin{equation}}
\newcommand{\en}{\end{equation}}
\newcommand{\eqm}{\begin{eqnarray}}
\newcommand{\enm}{\end{eqnarray}}
\newcommand{\eqmno}{\begin{eqnarray*}}
\newcommand{\enmno}{\end{eqnarray*}}
\newcommand{\eqml}[1]{\eql{#1}\begin{array}{rcl}}
\newcommand{\enml}{\end{array}\en}
\newcommand{\eql}{\begin{equation}\label}
\newcommand{\eqsub}[1]{\begin{subequations}\label{#1}\eqm }
\newcommand{\ensub}{\enm\end{subequations}}

\newcommand{\half}{{\frac{1}{2}}}

\newcommand{\goto}{\rightarrow}

\newcommand\lam{\lambda}
\newcommand\reals{{{\rm l} \kern -.15em {\rm R} }}

\newcommand\ddt{\frac{\partial}{\partial t}}
\newcommand\ddx{\frac{\partial}{\partial x}}
\newcommand\ddy{\frac{\partial}{\partial y}}
\newcommand\ddz{\frac{\partial}{\partial z}}

\newcommand\Dx{\Delta x}

\newcommand\Dt{\Delta t}

\newcommand\w{\tilde u}

\newcommand\Dq{\Delta \Q}

\newcommand\DtDx{\frac{\Delta t}{\Delta x}}
\newcommand\W{{\cal W}}

\newcommand\tW{\widetilde {\cal W}}

\newcounter{equationgroup}
\newcommand\Q{Q}

\renewcommand{\w}{w}

\newcommand{\tF}{\tilde F}

\renewcommand{\Dq}{\Delta\Q}

\newcommand{\imh}{_{i-1/2}}

\newcommand{\iph}{_{i+1/2}}

\newcommand{\ico}[1]{#1 \kern -1ex \raisebox{1.1ex}{$\circ$}}
\newcommand{\icos}[2]{#1 \kern -1.1ex \raisebox{1.1ex}[.5em]{$\circ$}^{#2}}

\newcommand{\inv}{^{-1}}

\newcommand{\Proj}{{{\rm l} \kern -.15em {\rm P} }}

}

\newcommand{\alert}[1]{{\bf \color{red} [#1]}}
\newenvironment{choice}{\left\{ \begin{array}{ll}}{\end{array}\right.}

\title{Computational Models of Material Interfaces for the Study of 
Extracorporeal Shock Wave Therapy}

\author{Kirsten Fagnan \thanks{Lawrence Berkeley National Laboratory, 1 Cyclotron Road, Berkeley, 
CA 94720} \and Randall J. LeVeque \thanks{Department of Applied Mathematics, University of 
Washington, Box 352420, Seattle, WA 98195} \and Thomas J. Matula \thanks{Applied Physics 
Laboratory, University of Washington, Seattle, WA 98195}}

\begin{document}

\maketitle

\begin{abstract}
Extracorporeal Shock Wave Therapy (ESWT) is a noninvasive treatment
for a variety of musculoskeletal ailments.
A shock wave is generated in water
and then focused using an acoustic lens or reflector so the energy of
the wave is concentrated in a small treatment region where mechanical stimulation
enhances healing.
In this work we have computationally
investigated shock wave propagation in ESWT by solving a Lagrangian
form of the isentropic Euler equations in the fluid and linear
elasticity in the bone using high-resolution finite volume methods.
We solve a full three-dimensional system of equations and use adaptive mesh
refinement to concentrate grid cells near the propagating shock.  
We can model complex bone
geometries, the reflection and mode conversion at interfaces, and the
the propagation of the resulting
shear stresses generated within the bone.   We discuss the validity of our
simplified model and present results validating this approach. 
\end{abstract}

\pagestyle{myheadings}
\thispagestyle{plain}

\section{Introduction}\label{sec:intro}
Extracorporeal shock wave therapy (ESWT) is a noninvasive treatment
for musculoskeletal conditions such as bone fractures
that fail to heal (non-unions), necrotic wounds, and strained
tendons\cite{cjwang_hip,eswt_shoulder}.  In this treatment a shock wave is generated in water and
then focused using an acoustic lens or reflector so that the energy
of the wave is concentrated in a small treatment region.  This
technique has been used since the 1980's, more widely in Europe and
Asia than in the US, where it is still considered experimental and
has limited FDA approval.  

Although the underlying biological
mechanisms are not well understood \cite{ogden}, the mechanical
compressional and/or shear
stress caused by the propagating shock wave is thought to
stimulate healing \cite{ogden,eswt_vegf,huang2010, morgan2008, turner1998, park1998,prendergast1997,lacroix2002,claes1999,isaksson2006,carter1998,goodship1985}.  There have been a variety of experimental studies
to investigate the effects of extracorporeal shock waves on urinary stones of varying hardness, as well as other hard tissues such as bone and its surrounding tissues (cartilage, tendon and fascia).  In particular, it has been shown that high extracorporeal shock wave energy actually fractures rat bones, but lower applied energy levels stimulated osteogenesis \cite{orthotripsy_ogden, valchanou, wolf2001}  We should also note that in a recent review article, Zelle, et. al. concluded that the majority of clinical studies were done at too high a level to enable a clinical recommendation \cite{zelle2010}.  However, they did conclude that shock wave therapy seems to stimulate the healing process in delayed unions and nonunions \cite{eswt_evidence_2010}. 
There are a number of other biological mechanisms in addition to stress that
potentially play a role in the body's response to ESWT. The focus of this
study, however, is on mechanical stress deposition and computational tools
for studying this aspect.

The medical shock wave devices
are similar to those used for extracorporeal shock
wave lithotripsy (ESWL), a widely-used non-surgical treatment for
kidney stones in which the focused shock waves have sufficient
amplitude to pulverize the kidney stone.  In shock wave therapy the
amplitudes are generally smaller and the goal is mechanical stimulation
rather than destruction, although in some applications such as the
treatment of heterotopic ossifications (HO) (see \Sec{ho}) larger
amplitudes may be used.

\begin{figure}
\begin{center}
a) \includegraphics[scale=0.35]{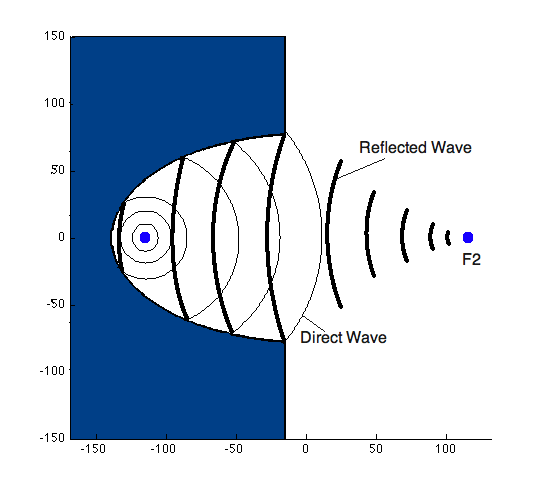}\hspace{1mm}
b) \includegraphics[scale=0.35]{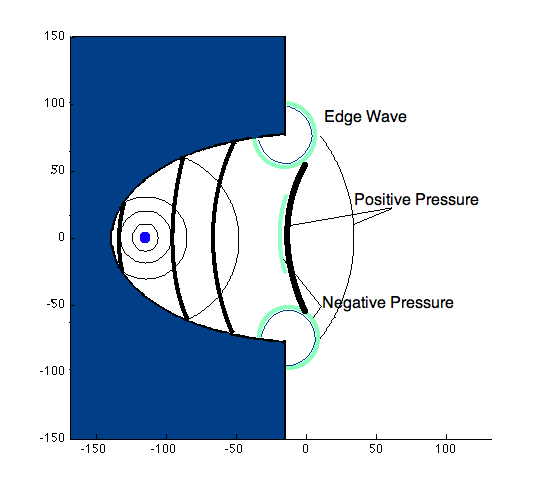}
\end{center}
\caption{Cartoon of the Dornier HM3 Lithotripter.  In a) the spherical wave is generated at F1, reflects off 
the ellipsoid and the reflected wave focuses at F2.  In b) the diagram illustrates the creation of the edge 
waves at the corner of the ellipsoid and the contribution of negative pressure to the tail of the ESWT 
pressure wave.}
\label{fig:hm3_diagram}
\end{figure}

Figure \ref{fig:hm3_diagram} shows the geometry of a laboratory
shock wave device 
modeled on the clinical Dornier HM3 lithotripter.  The three-dimensional
axisymmetric geometry consists of an ellipsoidal reflector made out
of metal and a cavity filled with water.  A spark plug at the focus
of the ellipse marked F1 generates a bubble which collapses and
creates a spherical shock wave that reflects and focuses at F2.  
The major and minor axes of the
ellipsoid in the HM3 are $a=140$mm and $b=79.8$mm, respectively.
The foci of this ellipse are at $(\pm 115, 0, 0)$ and the reflector
is truncated at $100$mm from F1, or $(-10,0,0)$.  

In the laboratory, this reflector is immersed in a bath of water
and objects can be
placed at the second focus of the ellipsoid, F2.
This device is in use at the Center for Industrial and Medical Ultrasound
(CIMU) at the University of Washington Applied Physics Laboratory and we
have used this geometry in order to compare directly with some laboratory
experiments.  Some preliminary comparisons were presented in 
\cite{kfagnan_hyp06}.

Computationally, we
use this geometry to calculate the initial condition by solving
two-dimensional axisymmetric Euler equations with the Tammann equation of
state (see \Sec{model_equations}).   These initial conditions are then fed into a full 
three-dimensional calculation near the focus at F2.

In addition to the HM3, we have also used the geometry of the
hand-held Sanywave device used in clinical studies
by our collaborator Dr. Michael Chang.
Some sample calculations related to the study of HOs are presented in
\Sec{ho}.

In each case, the ESWT pressure wave form
that is generated has a similar shape. There is a sharp increase
in pressure from atmospheric pressure ($\sim0.1$MPa)
to a peak  pressure ranging from 35 to 100
MPa over a very short rise time ($\sim10$ ns), 
followed by a decrease in pressure to $\sim-10$
MPa over $\sim5 \mu$s. The negative fluid pressure in the tail can lead to
cavitation bubbles, as discussed below.

Computational models for shock wave propagation and focusing can
aid in the study of ESWT.  
In particular, there are many open
questions concerning the interaction of shock waves with complex
three-dimensional geometries such as bone embedded in tissue.

Because of the difference in material properties, a wave hitting
the tissue/bone interface will be partially reflected, and the
transmitted wave will have a modified strength and direction of
propagation.  This can greatly affect the location and size of the
focal region as well as the peak pressure amplitude.  Moreover,
although the shock wave is primarily a pressure wave in soft tissue
(which has a very small shear modulus), at a bone interface 
mode conversion takes place and shear waves as well as compressional
waves are transmitted into the bone, generating a dynamically applied load.  

Bone healing is thought to be regulated in part by mechanical factors \cite{morgan2008, huang2010, turner1998, robling2002,huang2010,goodship1985,turner1998}.  Several studies have shown that the application of cyclic compressive and shear displacements can enhance healing through increased callus formation and ossification \cite{morgan2008, robling2002, roblingcyclic2002, saxon2005,park1998,weinbaum1994}.  The results also indicate that treatment is also dependent upon the rate, mode and magnitude of the stress deposition\cite{morgan2008}, as well as the gap size \cite{claes1995}. 

Carter, et. al. \cite{carter1998}, as well as, and Claes and Heigele \cite{claes1999}, proposed a model for skeletal tissue development based on hydrostatic pressure and tensile displacements \cite{claes1999}.  Other research has proposed a different model for skeletal tissue formation based on shear strain and fluid flow \cite{prendergast1997,lacroix2002}.  Augat, et. al. \cite{augat2001}, found that tensile displacements are not effective in enhancing bone formation.  This was further validated whe Isaksson, et. al. \cite{isaksson2006} investigated the models in \cite{claes1999, carter1998, prendergast1997, lacroix2002} and found that shear strain and fluid flow, were more accurate predictors of bone growth. However, no single model was able to predict certain features of the bone formation and healing process \cite{morgan2008}, highlighting the need for further research in this area.

The shear waves generated at the fluid/solid interface have also been shown to be important in 
the effective break up of kidney stones \cite{bailey_oleg,freund}.  An additional effect of ESWT is the
formation and collapse of cavitation bubbles that can cause tissue
damage.  While the shock wave is a compression wave, it is followed
by a rarefaction wave of expansion, and in the tail the fluid
pressure typically drops to negative values.  Reflection at interfaces
can lead to enhanced regions of expansion and to sufficiently
negative pressures that cavitation bubbles can form \cite{matula_direct, tanguay, coleman}.

To better understand all of these effects, it is desirable to have a
three-dimensional computational model that can simulate the focusing of
nonlinear shock waves and their interaction with arbitrarily complex
interfaces between different materials.  

\ignore{Many finite element models have been developed to estimate the mechanical stimuli computationally, where the user supplies the applied loads as well as the callus and bone geometry.  Our approach is similar to this, in that our applied load is the shock wave generated by the ESWT device and we are simulating the interaction with complex bone geometries.  By modeling the generation of the shock wave and the propagation through the bone and fluid, we are able to more easily simulate the applied load from ESWT with our bone geometries, something that is not as easily done in the other systems. }

In this paper we present an approach to this problem that has allowed
the study of some of these issues in a simplified context.  In
particular, we consider an idealized situation in which soft tissue
is replaced by water, ignoring its viscoelastic properties, and
modeled by the nonlinear compressible Euler equations with the
Tammann or Tait equation of state.  This has been used for prior ESWT 
work in water as well as biological-like materials \cite{ivings_toro, nakahara}.
Bone is modeled as an isotropic
and homogeneous linear elastic material \cite{fung,keaveny}.  

In reality, soft tissue and bone are very complex multiscale materials
with microstructures, inhomogeneities, and anisotropic properties.
Any attempt to model the biological effect of shock wave propagation
through such materials may require a more sophisticated and detailed
model than used here.  However, we believe that many of the macro-scale
shock propagation issues discussed above can be adequately and most
efficiently studied with a simplified model of the form considered here,
since the dominant effect we hope to capture is the reflection and
transmission of waves at interfaces between materials.  

The compressible Euler equations with the Tammann equation of state (see
\Sec{euler}) in two-dimensional axisymmetric geometry is used to model the
initial formation of the focusing shock wave.  These initial conditions are
then fed into a code that uses a simpler nonlinear
model, the Tait equation of state,
in a three-dimensional simulation of the fluid.  The compressible fluid
equations are written using a Lagrangian formulation that easily couples to
the isotropic linear elasticity equations used in the bone-like material.
The resulting equations have the same form everywhere, with a different
stress-strain relationship in the different materials.

A high-resolution finite volume method is used to solve these equations.
We use the wave-propagation algorithms described in \cite{rjl_book}
and implemented in Clawpack \cite{claw.org.url}.  These are Godunov-type
methods for the hyperbolic system that use solutions to the Riemann problem
between adjacent grid cells to determine a set of waves used to update the
solution, and second-order correction terms with slope limiters are added to
resolve the nearly discontinuous shock waves with minimal smearing or
nonphysical oscillation.

These methods are used on a purely rectangular Cartesian grid.  Each
grid cell has associated with it a set of material parameters
determining the material in the cell, in a unified manner so that
both fluid and solid can be modeled.   Complex geometry is handled
by using appropriate averaged values of these parameters in cells
that are cut by the interface.  This is described further in \Sec{interfaces}.
Averaging across the interface works quite well when the material properties
are sufficiently similar and in \Sec{interfaces} we show that this is the case even
for fluid/solid boundaries of the type we consider.  

We also use patch-based adaptive mesh refinement (AMR) to concentrate
grid cells in regions where they are most needed to resolve features
of interest.  The Clawpack software contains AMR software in both
two and three space dimensions and this software has been used
directly for the two-dimensional axisymmetric computations of the
initial shock wave described in \Sec{results}.  For the three-dimensional
problem we have used ChomboClaw \cite{chomboclaw}, an interface between
Clawpack and the Chombo code developed at LBL \cite{chombo}, which
provides an implementation of AMR on parallel machines using MPI.
Using ChomboClaw, the code originally developed using Clawpack was
easily converted into a code that was run on an NSF TeraGrid machine at
Texas Advanced Computing Center (TACC) and tested using up to 128 processors.

Extensive laboratory experiments have been performed on
shock wave devices to measure the wave form of shock waves produced by
various devices, the shape of the focal region, the peak amplitudes of
pressure observed in these regions, and other related quantities.  Most of
these experiments have been done in a water tank where the shock wave
propagates and focuses in a homogeneous medium where measurements are easily
done, or with phantoms (acrylic objects with well understood photoelastic properties)
 that are placed in the water as a proxy for bones or
kidney stones, with instrumentation such as pressure gauges or photographs
used to explore the interaction of the shock wave with the object.  In some
cases high-speed photographs of the shock wave have been obtained.
Creating phantoms from clear bi-refringent materials and using polarized
light it is even possible to photograph the shock wave propagating through
the object \cite{oleg_bailey}.  We have used some of these experiments to help
validate our numerical approach \cite{kfagnan_hyp06}.

Other researchers have also developed computational models for 
shock wave therapy and lithotripsy.
In prior work the pressure field has been modeled using linear and
nonlinear acoustics as well as the Euler equations with the Tait
equation of state.  Hamilton \cite{hamilton} used linear geometrical
acoustics, which holds under the assumption of weak shock strength,
to calculate the reflection of the spherical wave.  The diffraction
of the wave at the corner of the reflector was calculated using the
Kirchoff integral method.  Christopher's \cite{christopher_hm3}
model of the HM3 lithotripter used Hamilton's result as a starting
point and considered non planar sources.  Coleman et al \cite{coleman},
Averkiou and Cleveland \cite{cleveland_averkiou} used models based
on the KZK equation.  Tanguay \cite{tanguay} solved the full Euler
equations and incorporated cavitation effects as well as the edge
wave.

Our approach differs from these in that we consider the wave
propagation in both the fluid and solid by solving a single set of
equations that can model both materials.  This approach allows us
to investigate not only compression and tension effects of ESWT,
but also the propagation of shear waves in the solid.  Sapozhnikov
and Cleveland \cite{oleg_cleveland} have investigated the effect
of shear waves on spherical and cylindrical stones using linear
elasticity with a plane wave initial condition.  This initial
condition is an unfocused wave, which yields good results for small
objects, but would fail to capture the full ESWT pressure wave
interaction with three-dimensional bone geometries.

\section{Model equations}
\label{sec:model_equations}
To accurately model shock wave formation and propagation it is generally necessary to use nonlinear
equations of  compressible flow.  
In this work we use nonlinear equations for compressible liquids in the
fluid domain (water or soft tissue) and linear elasticity in the solid domain (bone).
The nonlinear compressible equations are written in a Lagrangian framework in terms of a reference
configuration, as is done for the linear elasticity equations.  This allows both sets of equations to be
written in the same form.  We apply finite volume methods to this form of the equations so that a single
computational grid (or set of nested grids with AMR) 
can be used over the entire domain.  Interfaces between fluid and
solid are represented by choosing averaged material parameters in each grid cell, as discussed further 
in \Sec{interfaces}.  

The system of equations we solve has the general form of a hyperbolic system 
of 9 equations
\begin{equation}\label{eqn:3dlgeulsys}
	q_t + f(q,x,y,z)_x + g(q,x,y,z)_y + h(q,x,y,z)_z = 0, 
\end{equation}
where the vector $q$ consists of the 6 components of the symmetric
strain tensor followed by the momenta,
 and the fluxes in general may be spatially varying based on material properties:
\begin{equation}\label{qfgh}
q = \bcm \epsilon^{11} \\ \epsilon^{22} \\ \epsilon^{33} \\ \epsilon^{12} \\ \epsilon^{23} \\ \epsilon^{13} \\  
\rho u\\ \rho v\\ \rho w\ecm, \quad
f(q,x,y,z) = \bcm u \\ 0 \\ 0 \\ v/2 \\ 0 \\ w/2 \\ \sigma^{11} \\ \sigma^{12} \\ \sigma^{13} \ecm, \quad
g(q,x,y,z) = \bcm 0\\ v \\ 0 \\ u/2 \\ w/2 \\ 0 \\ \sigma^{12} \\ \sigma^{22} \\ \sigma^{23} \ecm, \quad
h(q,x,y,z) = \bcm 0\\ 0\\ w\\ 0 \\ v/2 \\ u/2 \\ \sigma^{13} \\ \sigma^{23} \\ \sigma^{33} \ecm.
\end{equation}
In these expressions, $\rho= \rho(x,y,z)$ is the density of the material (the ``background density'' independent of the wave
propagating through the material) and the stress tensor $\sigma = \sigma(q,x,y,z)$ is in general a spatially varying 
function of $q$, linear in the solid and nonlinear in the fluid.  

\ignore{
\begin{align}
(\epsilon^{11})_t &= \frac{\partial u}{\partial x}\notag\\
(\epsilon^{22})_t &= \frac{\partial v}{\partial y}\notag\\
(\epsilon^{33})_t &= \frac{\partial w}{\partial z}\notag\\
(\epsilon^{12})_t &= \frac{1}{2}\left(\frac{\partial u}{\partial y} + \frac{\partial v}{\partial x}\right)\notag\\
(\epsilon^{23})_t &= \frac{1}{2}\left( \frac{\partial v}{\partial z} + \frac{\partial w}{\partial y}\right)\label{eqn:3dlgeulsys}\\
(\epsilon^{13})_t &= \frac{1}{2}\left( \frac{\partial u}{\partial z} + \frac{\partial w}{\partial x}\right)\notag\\
\rho u_t &= \frac{\partial \sigma^{11}}{\partial x} +\frac{\partial \sigma^{12}}{\partial y} +  \frac{\partial 
\sigma^{13}}{\partial z}\notag\\
\rho v_t &=   \frac{\partial \sigma^{12}}{\partial x} +\frac{\partial \sigma^{22}}{\partial y} + \frac{\partial 
\sigma^{23}}{\partial z}\notag \\
\rho w_t &=  \frac{\partial \sigma^{13}}{\partial x} + \frac{\partial \sigma^{23}}{\partial y} + \frac{\partial 
\sigma^{33}}{\partial z},\notag
\end{align}
where $\epsilon$ is the $3\times 3$ symmetric strain tensor and $\sigma(\epsilon)$ is the $3\times 3$
symmetric stress tensor.
}
Within the fluid domain $\sigma = -pI$, where $p$ is the scalar pressure and $I$ is
the identity matrix.  The pressure is a nonlinear function of the strain as discussed further below.
In the solid domain, $\sigma$ is a linear function of $\epsilon$ and is non-diagonal, allowing us to model
the propagation of shear waves as well as compressional waves.

In \Sec{euler} below we present the compressible fluid equations in their standard Eulerian form (the 
Euler
equations) and discuss two possible equations of state, the Tammann EOS and the simpler Tait EOS in 
which
the pressure is a function of density (or strain) alone, allowing us to drop the energy equation from the
Euler equations.   Then in \Sec{lagrangian} we rewrite these equations in the Lagrangian form
given above.  This can be done when modeling ESWT because the deformations are sufficiently small 
that 
the geometric nonlinearity of the equations can be ignored, adopting a Lagrangian frame and only 
considering
the nonlinearity of the stress-strain relation as given by the equation of state.  

In \Sec{elasticity} we discuss the linear elasticity model use to model bone.

\subsection{Compressible fluids in Eulerian form} \label{sec:euler}
Much of the previous work on ESWT has been centered around the use of the Euler equations with the 
Tait or Tammann equations of state.  These equations of state are typically used for modeling 
underwater 
explosions like the spark plug source of the lithotripter device \cite{hamilton, ivings_toro}.  In this section 
we discuss the full Euler equations and proceed to show why the Tait equation of state is sufficient for 
modeling ESWT.  Since this equation of state is a function only of the density, and can be rewritten as
a function of strain, we show in \Sec{lagrangian} how it 
can be modeled within the framework of elasticity, which enables us to model both the fluid 
and solid with the single system of equations given above. 

In three space dimensions the Euler equations take the form
\ignore{
\begin{equation}
\label{eulereqn}
	\frac{\partial}{\partial t} q + \frac{\partial}{\partial x} f(q) + \frac{\partial}{\partial y} g(q) 
 + \frac{\partial}{\partial z} h(q)= 0,
\end{equation}
with 
}
\begin{equation}
\label{eulereqn}
\ddt	\left[ \begin{array}{c} \rho \\ \rho u \\ \rho v \\ \rho w \\ E \end{array} \right  ]
+ \ddx \left[ \begin{array}{c} \rho u \\ \rho u^2 + p \\ \rho uv \\ \rho uw \\u(E+p) \end{array} \right]
	+ \ddy \left[ \begin{array}{c} \rho v \\ \rho uv \\ \rho u^2 + p \\ \rho vw \\v(E+p) \end{array} \right]
	+ \ddz \left[ \begin{array}{c} \rho v \\ \rho uw \\ \rho vw \\ \rho w^2+p \\w(E+p) \end{array} \right]
= 0.
\end{equation}

The total energy is $E 
= \rho e + \frac{1}{2}(u^2 + v^2 + w^2)$.

Several of the problems we investigated are axially symmetric and this enabled us to reduce the 
three-dimensional equations to a two-dimensional form.  If we first rewrite the equations in cylindrical 
coordinates $(r,\theta,z)$ and assume no variation and zero velocity in the $\theta$-direction, the 
system we obtain is reduced to two variables, $r$ and $z$.  The equations are
\begin{equation}
	\frac{\partial}{\partial t} \left[ \begin{array}{c} \rho \\ \rho u_r \\ \rho u_v \\ E \end{array} \right] +
	\frac{\partial}{\partial r} \left[ \begin{array}{c} \rho u_r \\ \rho u_r^2 + p \\ \rho u_r w_z \\ u_r(E+p) 
\end{array} \right] + 
	\frac{\partial}{\partial z} \left[ \begin{array}{c} \rho w_z \\ \rho u_r w_z\\ \rho w_z^2 + p \\ w_z(E+p) 
\end{array} \right ] = 
	\left[ \begin{array}{c} -(\rho u_r)/r \\ -(\rho u_r^2)/r \\ -(\rho u_r v_z) \\ u_r(E+p)/r \end{array} \right ],  
	\label{eqn:2daxisymeuler}
\end{equation}
where $u_r$ and $w_z$ denote the velocities in the $r$ and $z$ directions, respectively.  These 
equations are of the same form as the two-dimensional Euler equations, 
with the addition of geometric source terms that 
are a result of the variable transformation.  The source terms are never evaluate at $r=0$ since we are using a finite volume method where quantities are evaluated at cell-centers, that is, the smallest value of $r$ in a calculation is $\Delta x/2$.  We prefer to keep the equations in conservation form, so they can be efficiently solved using finite volume methods. 

In order to solve the system \eqn{eulereqn} or (\ref{eqn:2daxisymeuler}), we need to close the 
system with a relation between the pressure and conserved variables.  
The Tammann EOS \cite{ivings_toro} is applicable to a wide range of liquids, even with very strong 
shock
waves.  This equation of state has the form
\begin{equation}
	p = p(\rho, e) = (\gamma - 1)\rho e - \gamma p_{\infty},
	\label{eqn:tammaneos}
\end{equation}
where $p$, $\rho$ and $e$ are the pressure, density and specific internal energy, respectively, while $
\gamma$
and $p_\infty$ are constants depending on the fluid.  If $p_\infty=0$ this is the standard EOS for an ideal
gas, with $\gamma$ generally satisfying $1<\gamma < 5/3$, while for water $\gamma \approx 7.15$ and 
$p_\infty
\approx 300 \text{MPa}$.  For sufficiently weak shocks, this can be approximated
by the Tait equation of state,
\begin{equation}
	p = p(\rho) = B\left[ \left(\frac{\rho}{\rho_0}\right)^n - 1\right], 
	\label{eqn:taiteos}
\end{equation}
where $B$ is a pressure term that is a weak function of entropy, but is typically treated as a constant, 
and corresponds to $p_\infty$ from \eqn{eqn:tammaneos} while $n$ corresponds to $\gamma$.
In our work we take $B = 300 \text{MPa}$ and $n = 7.15$.

\begin{figure}
\begin{center}
a)\includegraphics[height=1.5in, width=2in]{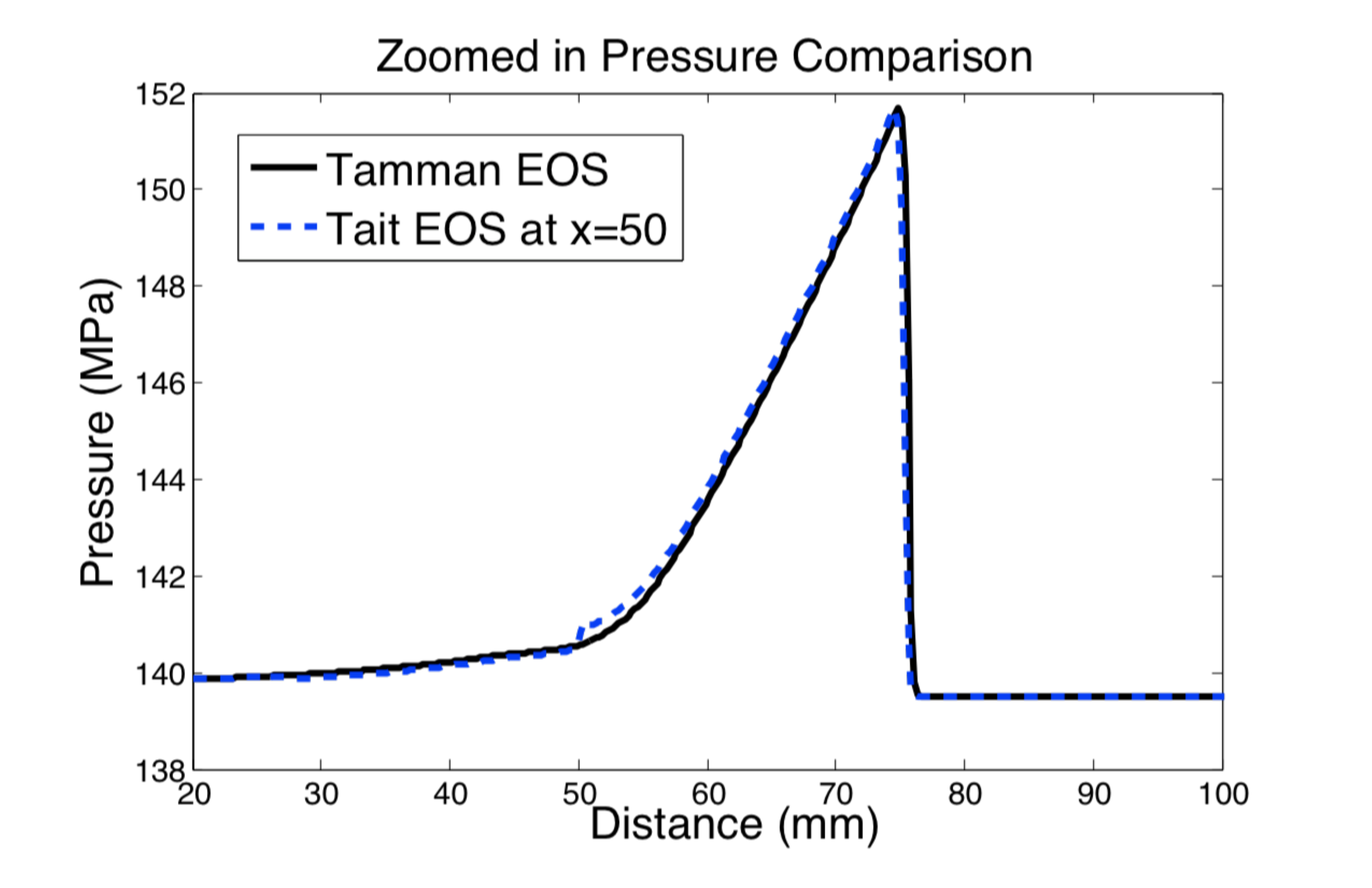}
b)\includegraphics[height=1.5in, width=2in]{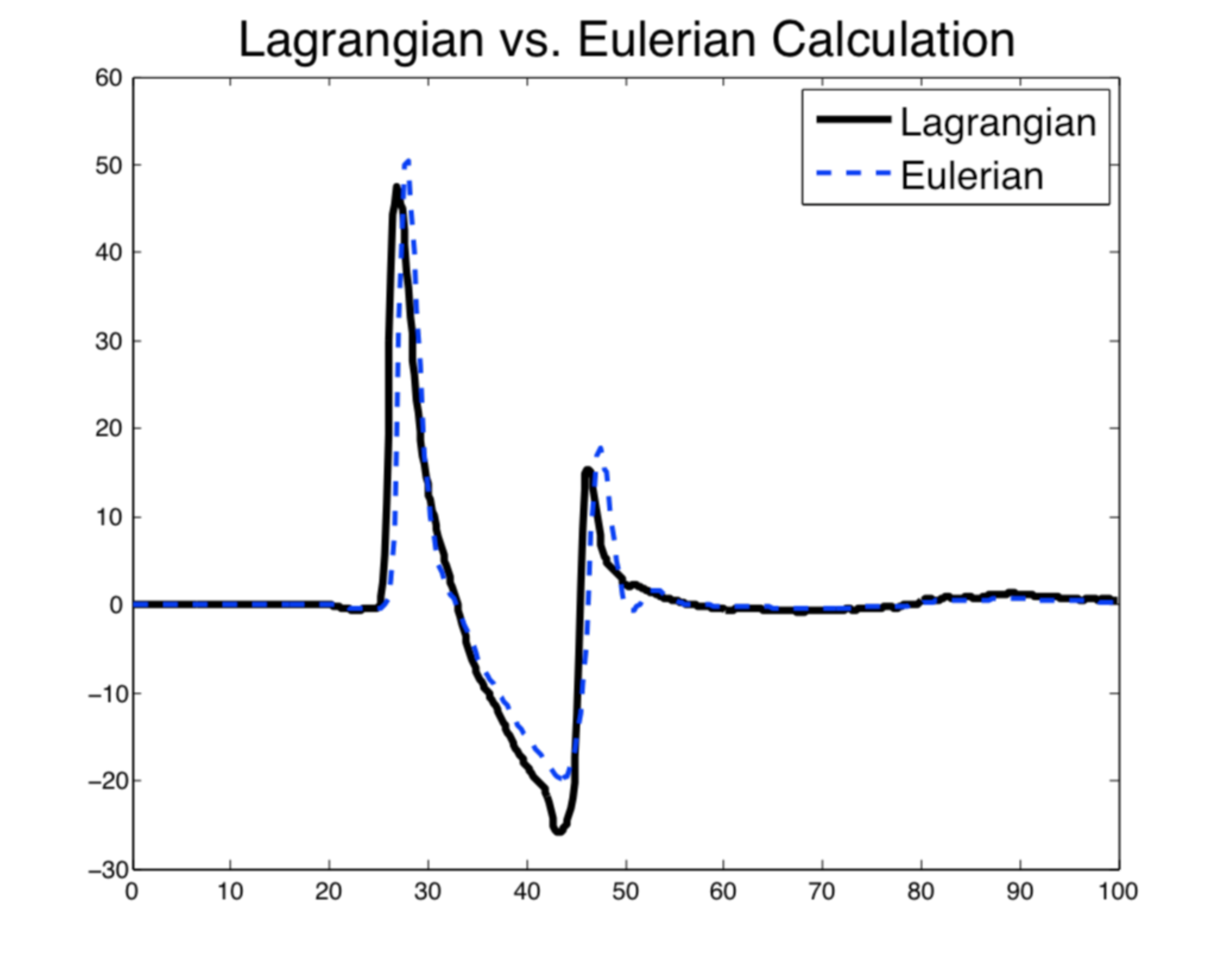}
\caption{a) Comparison of a pressure wave calculation performed using both the Tait (blue dashed curve) and Tammann (black curve) 
equations of state.  The results are nearly identical.  b) Comparison of the pressure pulse at F2 obtained 
in the Euler calculation (blue dashed curve) and the Lagrangian calculation (black curve).  It is clear that 
the two sets of equations give good agreement.  The wave in the Lagrangian case is slightly attenuated, 
but this may be due to error in initializing the calculation. In these calculations $\Delta x = 0.5 mm$.}
\label{fig:tait_tammann}
\end{center}
\end{figure}

It has been common practice to use the Tait EOS in shock wave therapy and lithotripsy 
models \cite{saito,nakahara}.  
This has been justified by noting studies that show that entropy changes across the shock are very small 
even
up to pressure jumps of 200 MPa \cite{nakahara}, which is beyond the range used in ESWT.
To verify this assumption, we performed computational experiments to compare the Tammann and Tait 
equations

of state for typical ESWT shock waves.  Since we have used the f-wave approach in our computational model, we can solve \ref{eqn:2daxisymeuler} with a spatially-varying equation of state.  We set up an experiment where the resulting shockwave (generated using the Tammann equation of state), was over 150 MPa.  Figure \ref{fig:tait_tammann} a) shows the results from this experiment.  The black solid curve is the result from solving with the Tammann EOS in the entire domain.  The blue dashed curve shows the result gotten by switching to the Tait EOS at $x=50$.  This enabled us to compare the two equations of state with the exact same initial condition.  There is a small disagreement at $x=50$ caused by a slight reflection at the interface due to the change in the equation of state.  Otherwise, the pressure profiles are nearly identical, giving
confidence that the calculations we are interested in can be done by solving the Euler equations with the 
Tait equation of state. This allows us to drop the equation for energy and obtain the simplified system
\begin{equation}
	\frac{\partial}{\partial t} \left[ \begin{array}{c} \rho \\ \rho u \\ \rho v \\ \rho w \end{array} \right] +
	\frac{\partial}{\partial x} \left[ \begin{array}{c} \rho u \\ \rho u^2 + p \\ \rho u v \\ \rho uw \end{array} \right] + 
	\frac{\partial}{\partial y} \left[ \begin{array}{c} \rho v \\ \rho u v\\ \rho v^2 + p \\ \rho vw \end{array} \right ] +
	\frac{\partial}{\partial z} \left[\begin{array}{c} \rho w \\ \rho uw \\ \rho vw \\ \rho w^2 + p \end{array} \right ] = 
	\left[ \begin{array}{c} 0 \\ 0 \\0 \\0 \end{array}\right ]. 
	\label{eqn:isentropic}
\end{equation}

\subsection{Compressible fluids in Lagrangian form} \label{sec:lagrangian}
In the case of a fluid where the shear modulus is zero, the stress tensor can be written as 
$\sigma(\epsilon) = -pI$, where $p$ is the pressure in the fluid and $I$ is the identity matrix.  
In the case of 
ESWT, the pressure only depends on changes in the density, and we can write $p(\epsilon)$ as a function 
of the strain tensor $\epsilon$.  Consider the movement of a material with
respect to a reference configuration and let $\delta = (\delta^x,~\delta^y,~\delta^z)$ be the
infinitesmal displacement.  

In three space dimensions, the full strain tensor is 
\begin{equation}
\epsilon = \left(\begin{array}{ccc} \delta^x_x & \frac{1}{2} (\delta^x_y + \delta^y_x) & \frac{1}{2} (\delta^z_y + \delta^y_z) \\
\frac{1}{2} (\delta^x_y + \delta^y_x) & \delta^y_y & \frac{1}{2} (\delta^x_z + \delta^y_x) \\
 \frac{1}{2} (\delta^z_y + \delta^y_z) & \frac{1}{2} (\delta^x_z + \delta^y_x) & \delta^z_z
 \end{array}
\right),
\end{equation}
where subscripts denote partial derivatives.

In the case of small deformations, we have that 
\begin{equation}
\rho = \frac{\rho_0}{1+tr(\epsilon)}
\end{equation}
where $\rho_0$ is the equilibrium density.

If we insert this into the Tait equation of state (\ref{eqn:taiteos}) we get
\begin{equation}\label{taiteps}
p(\epsilon) = B\left[ \left(\frac{1}{1 + tr(\epsilon)}\right)^n - 1\right].  
\end{equation}

Using the Lagrangian form is
only valid in the case where the displacements are small, so we calculated the maximum value of 
the displacements in a two-dimensional axisymmetric calculation with the Euler equations.  We found 
that for the maximum peak pressure of 50MPa, the corresponding maximum
velocity was $10^{-3}$m/s.  We then calculated the maximum displacement by 
integrating the velocity over the time of the calculation and found this to be on the order of 
$10^{-5}$mm.  The size of the grid cell is on the order of $10^{-1}$mm, so the displacements are
4 orders of magnitude smaller than the width of the grid cells.  
It is therefore reasonable to assume that the density in each grid cell is essentially 
constant and that the Lagrangian framework of the elasticity equations will be valid for 
the fluid.  

To test this, we took the same initial condition for the 2D axisymmetric Euler equations with the 
Tammann equation of state and the corresponding 2D axisymmetric Lagrangian form of the equations 
with the Tait equation of state and measured the pressure at the focus, F2.  The results in Figure 
\ref{fig:tait_tammann} b) demonstrate reasonably good agreement between the two cases, but the 
Lagrangian form is slightly attenuated.  This may be due to conversion of the initial condition from the 
conserved variables in the Euler equations (\ref{eqn:2daxisymeuler}) to those in the elasticity equations 
(\ref{eqn:3dlgeulsys}).

\begin{figure}
a)\includegraphics[width=2in]{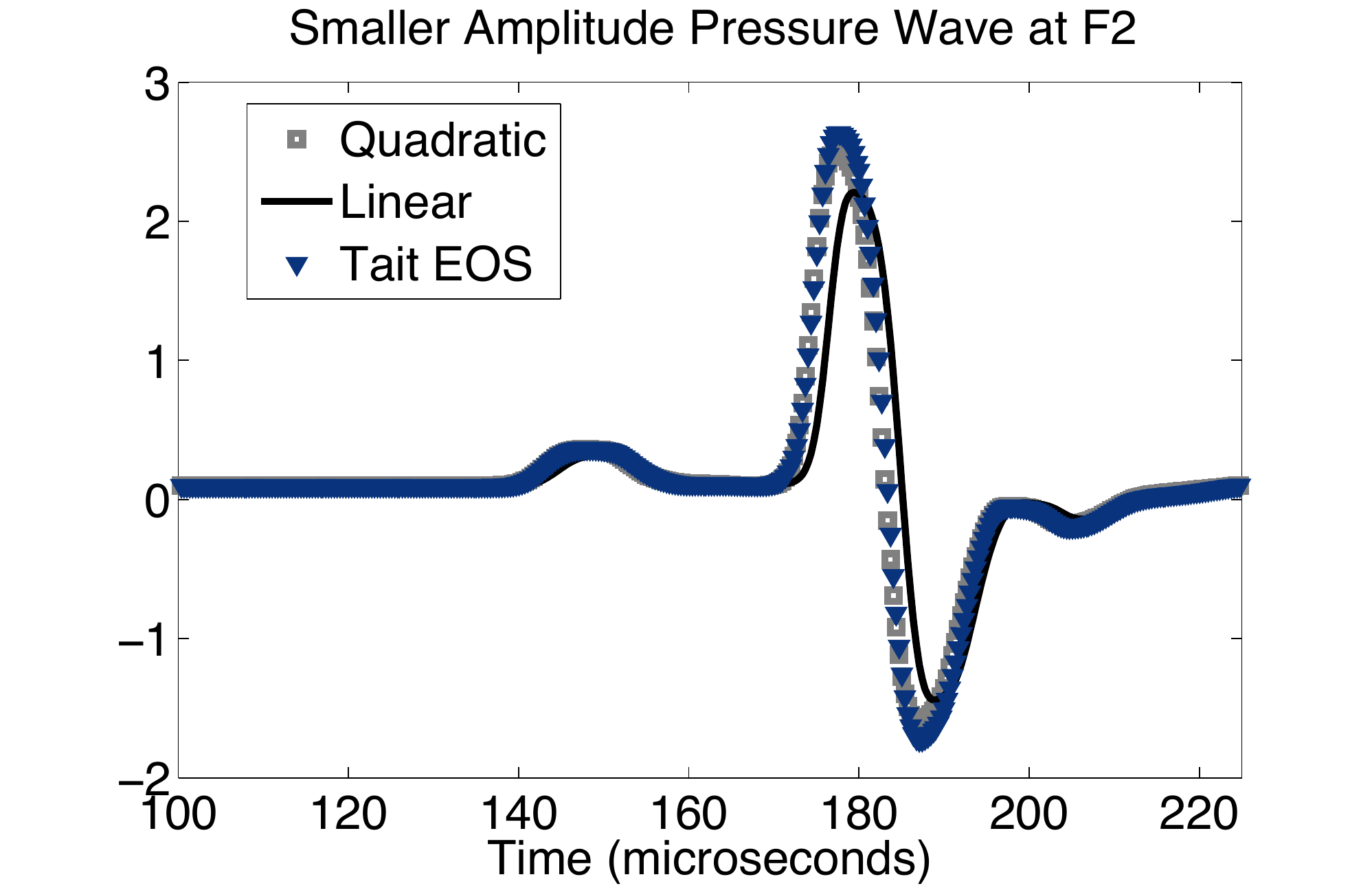}\hspace{-6mm}
b)\includegraphics[width=2in]{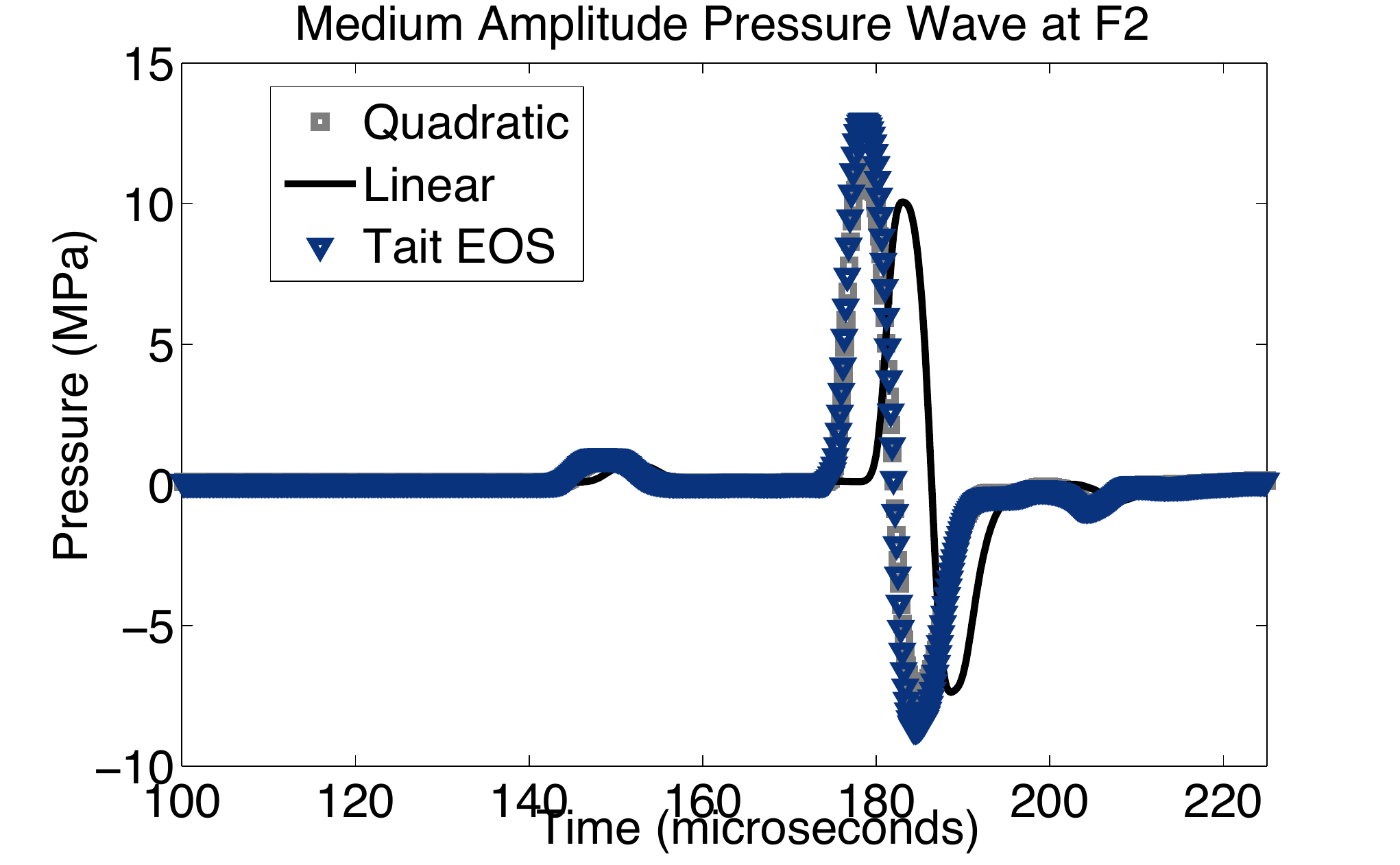}\hspace{-6mm}
c)\includegraphics[width=2in]{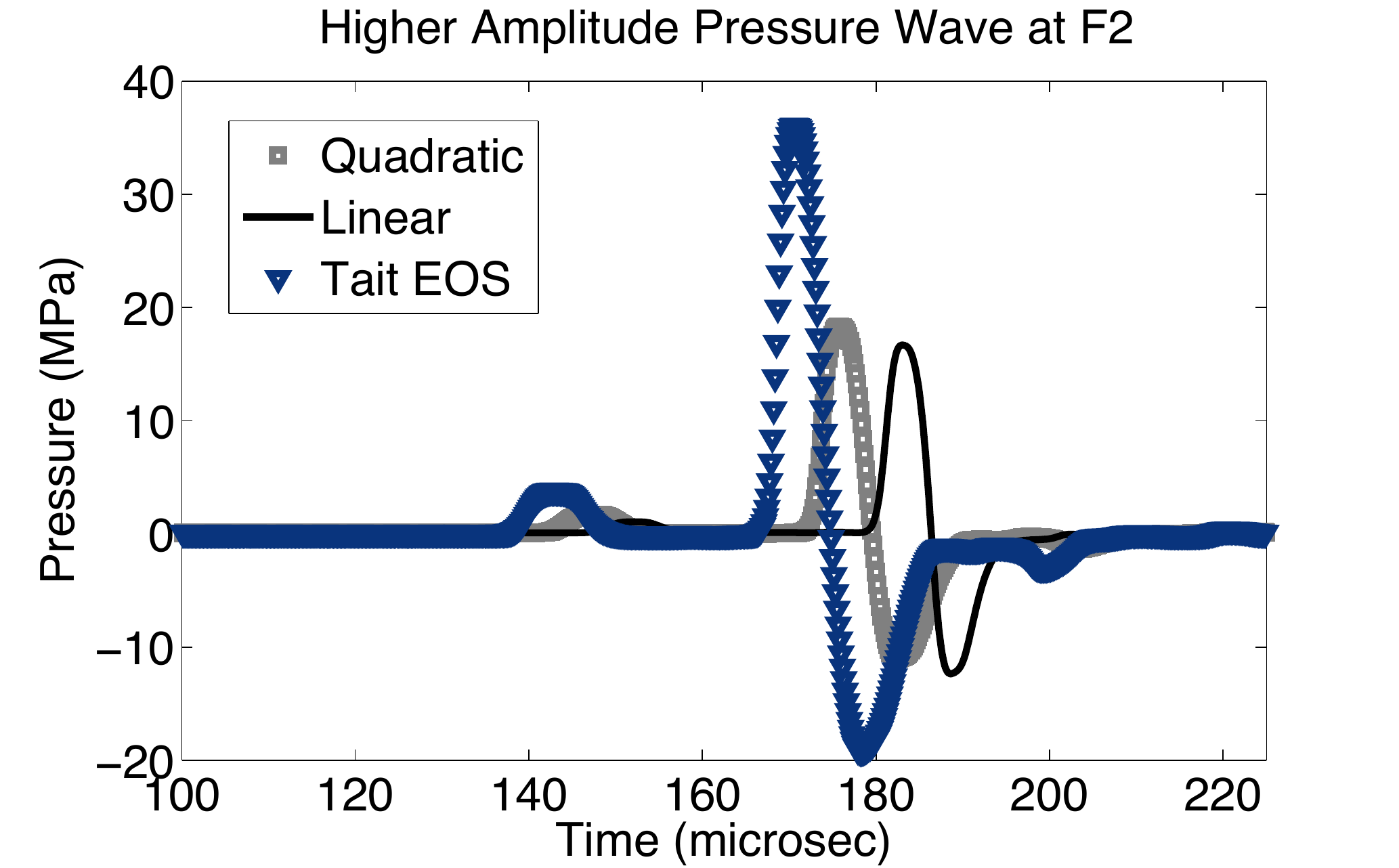}\hspace{-6mm}
\caption{This series of figures shows a pressure gauge measurement at F2 of different versions of the Tait EOS at different amplitudes.  The triangular markers indicate the full nonlinear Tait EOS, the solid line is a linearized version and the square markers are a quadratic version.  The linearized versions of the EOS work reasonably well at small amplitudes, but it's clear from figure c) that as the pressures increase to those observed in ESWT, the full nonlinear equation of state must be used.}
\label{fig:tait_linearization_test}
\end{figure}

Since the displacements are small, we also considered the possibility that nonlinearity in the fluid
could be ignored, so we could instead using a linearized version of the Tait 
equation of state.   Then we would be able to simply use the linear elasticity 
equations throughout the domain, in both the fluid and solid materials.  
If we assume a small perturbation to the strain, $\epsilon + \delta \epsilon$, we can expand the Tait EOS (\ref{eqn:taiteos}) as a Taylor series about $\epsilon$, 
\begin{equation}
	p(\epsilon + \delta \epsilon) = p_0 + p'(\delta \epsilon) \epsilon + \frac{p''(\delta \epsilon)}{2} \epsilon^2 +
\cdots
	\label{eqn:tait_taylor}
\end{equation}
If we keep the first two terms of the expansion, the EOS has been linearized and we will call this the 
``linear Tait EOS.''  Similarly, we will refer to the equation obtained by keeping the first three terms of the expansion as the ``quadratic Tait EOS.''
One-dimensional tests of both possibilities are shown in 
Figure~\ref{fig:tait_linearization_test}, for three different wave amplitudes.
For a wave with maximum amplitude less than 3 MPa there is fairly good agreement, however, as the amplitude is increased, as is required for ESWT, the linear and quadratic equations of state do not capture the correct behavior.  Thus we used the full Tait EOS in the fluid domain.
\ignore{
We found that for smaller amplitude waves, the linear Tait EOS is insufficient 
for predicting the correct arrival time or amplitude.  If we include the quadratic term, then we can 
accurately determine the arrival time, but the amplitude is still incorrect.  As we increase the amplitude of 
the pressure wave, i.e. there is a larger jump in pressure across the shock, neither the linear nor 
quadratic EOS can capture the correct arrival time or amplitude.  This is demonstrated in Figure
\ref{fig:tait_linearization}.
}

\ignore{
\subsection{Comparison of Lagrangian and Eulerian Solution}
\label{sec:lag_eul_comp}
Given the result in Figure \ref{fig:tait_tammann} a), combined with the small displacements in our 
model, we are confident that the two systems of equations should agree for the weak shocks present in 
ESWT.  We performed an axisymmetric calculation with the Euler equations using the Tammann 
equation of state to generate an initial condition for the experiment.  We then measured the focused 
wave form at $x=115.0$ for both the Eulerian and Lagrangian form of the equations.  In Figure 
\ref{fig:tait_tammann} b) we see that the two models give nearly the same solution, although the wave 
from the Lagrangian calculation has been slightly attenuated.
}

\subsection{Elasticity equations}
\label{sec:elasticity}
\ignore{
In the case of infinitesimal or small displacements, it is possible to write the elasticity equations as 
\begin{align}
 \epsilon_t - u_x &= 0 \label{eqn:nlelast}\\
 \rho_0 u_t - \sigma(\epsilon)_x &= 0  \notag\\
\end{align}
where the higher order terms in the deformation tensor have been neglected.  This form of the equations 
still allows for nonlinear behavior to be incorporated through the stress-strain relationship, $
\sigma(\epsilon)$.  We take advantage of this formulation in order to incorporate the nonlinear Tait 
equation of state in the fluid. 

In the case of a fluid where the shear modulus is zero, the stress tensor can be written as $
\sigma(\epsilon) = -pI$, where $p$ is the pressure in the fluid and $I$ is the identity matrix.  In the case of 
ESWT, the pressure only depends on changes in the density, and we can rewrite $p(\rho)$ as a function 
of the strain tensor $\epsilon$.  In the case of small deformations, we have that 
\begin{equation}
\rho = \frac{\rho_0}{\det(F)} = \frac{\rho_0}{1+tr(\epsilon)},
\end{equation}
where $F$ is the deformation gradient tensor.

If we plug this into the Tait equation of state (\ref{eqn:taiteos}) we get
\begin{equation}
p(\rho) = p(\epsilon) = B\left[ \left(\frac{1}{1 + tr(\epsilon)}\right)^n - 1\right].  
\end{equation} 
Therefore, in the fluid, we get the following for our stress-strain relationship 
\begin{equation}
\sigma(\epsilon) = -(B\left[ \left(\frac{1}{1 + tr(\epsilon)}\right)^n - 1\right]) I
\end{equation} 
}
In the current work we model bone as a 
linear isotropic solid.  We use the equations \eqn{eqn:3dlgeulsys}
together with Hooke's law
\begin{align}\label{hookes}
\sigma^{11} &= C_{11} \epsilon^{11} + C_{12} \epsilon^{22} + C_{13} \epsilon^{33} \\
\sigma^{22} &= C_{21} \epsilon^{11} + C_{22} \epsilon^{22} + C_{23} \epsilon^{33}  \\
\sigma^{33} &=C_{31} \epsilon^{11} + C_{32} \epsilon^{22} + C_{33} \epsilon^{33}  \\
\sigma^{12} &= C_{44} \epsilon^{12} \\
\sigma^{13} &= C_{55} \epsilon^{13}\\
\sigma^{23} &=  C_{66} \epsilon^{23}. \\
\end{align}
where the spatially-varying scalar coefficients $C_{ij}(x,y,z)$ are determined by the properties of the material being modeled.  
The parameters used for the bone model were found in \cite{martin_burr_sharkey}.  

For an isotropic material
we can relate the $C_{ij}$ above to the two Lam\'{e} parameters,
$\lambda$ and $\mu$, that are used to model different elastic
materials.  $C_{ii} = \lambda + 2 \mu$ for $i=1,\dots,3$, $C_{ii}=2
\mu$ for $i=4,\dots,6$, and $C_{ij} = \lambda$ for $i \ne j$.  Here
$\mu$ is the shear modulus and $\lambda + 2\mu$ is the bulk modulus
of the material.  Note that the $\lambda$ here is different from
the $\lambda^{i}$ used to denote the eigenvalues elsewhere in the
paper.

Linear elasticity has been used extensively in the literature to model both trabecular and cortical bone \cite{keaveny,fung}.  Linear viscoelastic models have also been used for ultrasound wave propagation in bone \cite{visc_bone}.
Our model could be extended to orthotropic models, requiring 9 material parameters, as has also been
used for bone modeling, e.g. \cite{orthotropic}.

\section{Eigenstructure of the hyperbolic system}
\label{sec:eigen}
The full three-dimensional system of equations \eqn{eqn:3dlgeulsys}
models both the nonlinear fluid and the linear elastic bone as described in the preceeding sections.
This system can be written in quasi-linear form
\ignore{
\begin{align}
(\epsilon^{11})_t &= \frac{\partial u}{\partial x}\notag\\
(\epsilon^{22})_t &= \frac{\partial v}{\partial y}\notag\\
(\epsilon^{33})_t &= \frac{\partial w}{\partial z}\notag\\
(\epsilon^{12})_t &= \frac{1}{2}\left(\frac{\partial u}{\partial y} + \frac{\partial v}{\partial x}\right)\notag\\
(\epsilon^{23})_t &= \frac{1}{2}\left( \frac{\partial v}{\partial z} + \frac{\partial w}{\partial y}\right)\label{eqn:
3dlgeulsys}\\
(\epsilon^{13})_t &= \frac{1}{2}\left( \frac{\partial u}{\partial z} + \frac{\partial w}{\partial x}\right)\notag\\
\rho u_t &= \frac{\partial \sigma^{11}}{\partial x} +\frac{\partial \sigma^{12}}{\partial y} +  \frac{\partial 
\sigma^{13}}{\partial z}\notag\\
\rho v_t &=   \frac{\partial \sigma^{12}}{\partial x} +\frac{\partial \sigma^{22}}{\partial y} + \frac{\partial 
\sigma^{23}}{\partial z}\notag \\
\rho w_t &=  \frac{\partial \sigma^{13}}{\partial x} + \frac{\partial \sigma^{23}}{\partial y} + \frac{\partial 
\sigma^{33}}{\partial z}.\notag
\end{align}
}

\begin{equation}\label{quasilin}
	q_t + A(q,x,y,z) q_x + B(q,x,y,z) q_y + C(q,x,y,z) q_z = 0, 
\end{equation}
where $A$,$B$ and $C$ are the Jacobians of the flux functions in the $x$, $y$ and $z$ directions 
respectively.  For the multi-dimensional methods implemented in Clawpack, we need the solution to 
the Riemann problem along slices in each coordinate direction.  Here we provide the details for the 
solution in the $x$-direction, but the solution in the $y$ and $z$ directions are similar with appropriate 
perturbations to the $B$ and $C$ matrices. 

The corresponding Jacobian for this system in the $x$-direction is
\begin{equation}
A(q,x,y,z) = \frac{\partial f(q,x,y,z)}{\partial x}
= -\left(\begin{array}{ccccccccc}0&0&0&0&0&0&\frac{1}{\rho_0}&0&0\\ 
0&0&0&0&0&0&0&0&0\\0&0&0&0&0&0&0&0&0\\0&0&0&0&0&0&0&\frac{1}{2\rho_0}&0\\
0&0&0&0&0&0&0&0&0\\0&0&0&0&0&0&0&0&\frac{1}{2\rho_0}\\ \sigma^{11}_{\epsilon^{11}} & 
\sigma^{11}_{\epsilon^{22}} &\sigma^{11}_{\epsilon^{33}} & 0 &0 &0 &0&0&0\\0&0&0&
\sigma^{12}_{\epsilon^{12}}&0&0&0&0&0\\0&0&0&0&0&\sigma^{13}_{\epsilon^{13}}&0&0&0 \end{array} 
\right)
\label{eqn:jac3dlgeul}
\end{equation}
where $\sigma^{11}_{\epsilon^{33}}$, for example, denotes the partial derivative of
$\sigma^{11}$ with respect to $\epsilon^{33}$.  In the linear elastic case this is
simply the coefficient $C_{13}$, but the above form also applies to the nonlinear
compressible equations.
The  spatial variation in  $f(q,x,y,z)$ and the Jacobian $A$ result from 
allowing the material parameters such as density and elastic moduli to vary
in space.
The Jacobians in the $y$- and $z$- directions are similar with the entries perturbed appropriately.

The eigenvalues for system (\ref{eqn:jac3dlgeul}) are
\begin{equation}
\lambda^{1,2} = \pm \sqrt{\frac{\sigma^{11}_{\epsilon^{11}}}{\rho_0}}; \quad \lambda^{3,4} = \pm 
\sqrt{\frac{\sigma^{12}_{\epsilon^{12}}}{2\rho_0}}; \quad\lambda^{5,6} = \pm 
\sqrt{\frac{\sigma^{13}_{\epsilon^{13}}}{2\rho_0}}; \quad\lambda^{7,8,9}=0.
\label{eqn:3dlgeuleval}
\end{equation}
When modeling a fluid where the shear stress is zero, there are seven
zero-speed eigenvalues since $\sigma^{12}_{\epsilon^{12}} = 
\sigma^{13}_{\epsilon^{13}} = 0$.  
Only the compressional waves corresponding to $\lambda^{1,2}$ propagate
with nonzero speed.  
Note that with the Tait equation of state \eqn{taiteps},
\begin{align}
\sigma^{11}_{\epsilon^{11}} =    \frac{\partial\sigma^{11}}{\partial\epsilon^{11}} &= B n \left( \frac{1}{(1+\epsilon^{11} + \epsilon^{22} + \epsilon^{33})} 
\right )^{n +1} \label{eqn:sigepstait}\\
			&=  \frac{n(p+B)}{(1+\text{tr}(\epsilon))}.\notag
\end{align}

In the small amplitude acoustic limit $\epsilon\goto 0$,
from (\ref{eqn:3dlgeuleval}) we obtain the wave speeds 
\begin{equation}
	\pm \sqrt{\frac{n (p+B)}{\rho_0}}.
\label{eqn:lagisenspeed}
\end{equation}
which are the expected waves speeds for compressional waves in the Lagrangian form with this
equation of state.

For the elastic solid, on the other hand, waves 1 and 2 correspond to P-waves while waves 4--6
correspond to S-waves, and the expected wave speeds are recovered based on the elastic
coefficients given in \Sec{elasticity}.  For example, in the $x$-direction the P-wave speeds are
\begin{equation}
	\pm \sqrt{\frac{C_{11}}{\rho_0}},
\end{equation}
and the S-wave speeds are
\begin{equation}
	\pm \sqrt{\frac{C_{44}}{2\rho_0}} ~~\text{  and  }~~ \pm \sqrt{\frac{C_{55}}{2\rho_0}}.
\end{equation}

The corresponding eigenvectors for system (\ref{eqn:jac3dlgeul}) are
\begin{equation}
\label{r1-6}
r^{1,2} = \left ( \begin{array}{c} 1\\0\\0\\0\\0\\0\\ \pm\sqrt{\rho_0 \sigma^{11}_{\epsilon_{11}}}\\ 0\\0
\end{array}
\right ),\qquad
r^{3,4} = \left ( \begin{array}{c} 0\\0\\0\\1\\0\\0\\0\\\pm \sqrt{2 \rho_0 \sigma^{12}_{\epsilon^{12}}}\\0
\end{array}
\right),\qquad
r^{5,6} = \left ( \begin{array}{c} 0\\0\\0\\0\\0\\1\\0\\0\\ \pm\sqrt{2 \rho_0 \sigma^{13}_{\epsilon_{13}}}
\end{array}
\right )
\end{equation}
for the P-waves and S-waves, and
\eql{r789}
r^7 = \left ( \begin{array}{c}-\sigma^{11}_{\epsilon^{22}} \\ {\sigma^{11}_{\epsilon^{11}}}\\0\\0\\0\\0\\0\\0\\0
\end{array}
\right ),\qquad
r^8 = \left ( \begin{array}{c}-\sigma^{11}_{\epsilon^{33}}\\0\\ {\sigma^{11}_{\epsilon^{11}}}\\0\\0\\0\\0\\0\\0
\end{array}
\right ),\qquad
r^9 = \left ( \begin{array}{c}0\\0\\0\\0\\1\\0\\0\\0\\0
\end{array}
\right )
\end{equation}
for the stationary waves.

\subsection{Axisymmetric form of the equations}
\label{sec:axisym}
We used the two-dimensional axisymmetric form of the equations to generate an initial condition for our 
three-dimensional calculations, as well as for validation of our model.  

The three-dimensional equations in cylindrical coordinates are
\begin{align}
\epsilon^{rr}_t &= \frac{\partial u}{\partial r}\notag\\
\epsilon^{\theta \theta}_t &= \frac{u}{r} + \frac{1}{r} \frac{\partial v}{\partial \theta}\notag\\
\epsilon^{zz}_t &= \frac{\partial w}{\partial z}\notag\\
\epsilon^{rz}_t &= \frac{1}{2} \left(\frac{\partial u}{\partial z} + \frac{\partial w}{\partial r}\right)\notag\\
\epsilon^{r \theta}_t &= \frac{1}{2} \left(\frac{\partial v}{\partial r} + \frac{1}{r} \frac{\partial u}{\partial \theta} - 
\frac{v}{r}\right)\notag\\
\epsilon^{\theta z}_t &= \frac{1}{2r}\left( \frac{\partial w}{\partial \theta} + \frac{\partial v}{\partial z}\right)\label{eqn:3dcyllageul}\\
\rho u_t &= \frac{1}{r} \frac{\partial \sigma^{r \theta}}{\partial \theta} + \frac{\partial \sigma^{rr}}{\partial r} + 
\frac{\sigma^{rr} - \sigma^{\theta \theta}}{r} + \frac{\partial \sigma^{rz}}{\partial z}\notag\\
\rho v_t &= \frac{1}{r} \frac{\partial \sigma^{\theta \theta}}{\partial \theta} + \frac{\partial \sigma^{r \theta}}
{\partial r} + \frac{2 \sigma^{r \theta}}{r} + \frac{\partial \sigma^{z \theta}}{\partial z}\notag \\
\rho w_t &= \frac{1}{r} \frac{\partial \sigma^{z \theta}}{\partial \theta} + \frac{\partial \sigma^{zz}}{\partial z} 
+ \frac{\partial \sigma^{rz}}{\partial r} + \frac{\sigma^{rz}}{r}.\notag
\end{align}

If we assume that $v = \epsilon_{\theta z} = \epsilon_{r \theta} = 0$ and there is no variation in the $\theta
$ direction, then the system (\ref{eqn:3dcyllageul}) simplifies to
\begin{align}
\epsilon^{rr}_t &= \frac{\partial u}{\partial r}\notag\\
\epsilon^{\theta \theta}_t &= \frac{u}{r} \notag\\
\epsilon^{zz}_t &= \frac{\partial w}{\partial z}\notag\\
\epsilon^{rz}_t &= \frac{1}{2}\left(\frac{\partial u}{\partial z} + \frac{\partial w}{\partial r}\right)\label{eqn:3daxisym}\\
\rho u_t &= \frac{\partial \sigma^{rr}}{\partial r} + \frac{\sigma^{rr} - \sigma^{\theta \theta}}{r} + \frac{\partial 
\sigma^{rz}}{\partial z}\notag\\
\rho w_t &= \frac{\partial \sigma^{zz}}{\partial z} + \frac{\partial \sigma^{rz}}{\partial r} + \frac{\sigma^{rz}}
{r}.\notag
\end{align}
It is interesting to note here that the strain in the $\theta \theta$ direction is non-zero and in this case is 
called the hoop strain.  A uniform radial displacement is not a rigid body motion, as it would be in the 
two-dimensional plane strain case, but instead produces a circumferential strain.  This is because the 
original circumference of the cylinder is $2\pi r$, but when there is a strain in the radial direction the 
circumference grows to $2\pi (r + u_r)$, inducing a strain $2\pi u_r / 2 \pi r = u_r/r$.

The Jacobian for system (\ref{eqn:3daxisym}) in the z-direction is
\begin{equation}
f^{\prime}(q) = -\left(
\begin{array}{cccccc} 0&0&0&0&\frac{1}{\rho_0}&0 \\ 0&0&0&0&0&0 \\ 0&0&0&0&0&0 \\ 0&0&0&0&0&\frac{1}{2\rho_0} \\ 
\sigma^{rr}_{\epsilon_{rr}}&\sigma^{rr}_{\epsilon_{zz}}&0&0&0&0\\ 0&0&0&
\sigma^{rz}_{\epsilon_{rz}}&0&0
\end{array}
\right), 
\end{equation}
and has an eigen-structure that is equivalent to the two-dimensional elasticity equations, with the 
addition of a second zero-speed eigenvalue.  

These equations have the structure
\begin{equation}\label{2dsrc}
q_t + f(q)_r + g(q)_z = S(q,r),
\end{equation}
with source terms
\begin{align}
\epsilon^{\theta \theta}_t &= \frac{u}{r}  \notag\\
\rho u_t &= \frac{\sigma_{rr} - \sigma_{\theta \theta}}{r}\label{eqn:axisymeulsource}\\
\rho w_t &= \frac{\sigma_{rz}}{r}.\notag
\end{align} 
In Clawpack, we solve these equations with a fractional-step method.  The full problem is split into two subproblems that are solved independently.  
We first solve the homogeneous system obtained by setting $S\equiv 0$ in \eqn{2dsrc} 
using the wave propagation algorithm described in \Sec{numerics}, and then solve
\begin{equation}
q_t = S(q,r),
\end{equation}
with an appropriate ODE solver.  For (\ref{eqn:axisymeulsource}), we use forward Euler. 

\section{Numerical methodology}
\label{sec:numerics}
We used the wave-propagation algorithms described in \cite{rjl_book}
and implemented in Clawpack \cite{clawpack} to solve the hyperbolic systems
of PDEs derived in the preceeding sections.
In this section we provide the basic details of the
numerical methodology and the approximate solution to the Riemann problem
with a spatially-varying flux function, similar to what was done in
\cite{rjl_nonlinear}  We also discuss computational issues that require
the use of adaptive mesh refinement.

\subsection{Riemann solvers and wave-propagation algorithms}
\label{sec:riemann_solver}
Recall that the ``Riemann problem'' 
is the initial value problem for a 1-dimensional hyperbolic system 
of the form
\eql{qtfqx}
	q_t + f(q,x)_x = 0, 
\end{equation} 
with special initial data consisting of two constant states separated by a discontinuity
\eql{Rpdata}
	    q_{0}(x) = \begin{choice} Q_l &\text{if}~x<0 \\
                    Q_r &\text{if}~x>0 \end{choice}.
\end{equation} 
If the flux function is spatially varying then we also use a piecewise-defined 
flux function with
\eql{Rpf}
	    f(q,x) = \begin{choice} f_l(q) &\text{if}~x<0 \\
                    f_r(q) &\text{if}~x>0 \end{choice}.
\end{equation} 
The Riemann problem plays a fundamental role in the theory and computation
of hyperbolic problems, since the Riemann solution consists of waves
propagating at constant speeds and can generally be computed.
For nonlinear systems of equations this is
often replaced by an approximate Riemann
solver as will be discussed below. 

For a linear system of equations
$q_t + A(x)q_x=0$ the Riemann solution is easily computed in terms of
the eigenvectors and eigenvalues of the matrices $A_l$ to the left of the
interface and $A_r$ to the right of the interface.
We begin by discussing the linear case with a constant matrix $A$ and turn
to the variable-coefficient (heterogeneous media) case in \Sec{hetero}.
We assume the matrix $A$ is diagonalizable,
\eql{Adecomp}
A = R\Lambda R\inv
\end{equation} 
where $R$ is the matrix of eigenvectors and $\Lambda$ is the diagonal matrix
of eigenvalues.  The Riemann solution is computed by decomposing
$\Delta Q = Q_r - Q_l$ as a linear combination of eigenvectors of $A$,
\eql{dqwaves}
\Dq = \sum_{p=1}^m \alpha^p r^p, \qquad\text{where}~~ \alpha = R\inv \Dq.
\end{equation} 
We denote the $p$th wave by $\W_p = \alpha^p r^p$, where $p=1,~2,~\ldots,~m$
and the number of waves $m$ is equal to the number of equations in the
system.

We use finite volume methods in which $Q_i^n$ represents a cell average of the
vector $q$ in cell $i$ at time $t_n$ (still in one space dimension).
In Godunov's method the cell average is updated by the waves entering the
cell from the interfaces to the left and the right, and each wave updates
the cell average by $\W^p$, the jump in $q$ across the wave, multiplied by the
distance the wave propagates over the time step and divided by the length of
the cell, i.e., the cell average is updated by $\frac{\lam^p \Dt}{\Dx}\W^p$.
To express the total update to a cell, it is convenient to 
define matrices $A^+$ and $A^-$ via
\eql{Apm}
A^\pm = R\Lambda^\pm R\inv, \qquad\text{where}~~\Lambda^\pm =
\text{diag}(\lambda_p^\pm)
\end{equation}
with $\lambda^+ = \max(\lambda, 0)$ and $\lambda^- = \min(\lambda, 0)$.
Then the cell average is updated by
\eql{wp1d1}
Q_i^{n+1} = Q_i^n - \DtDx (A^+\Dq\imh + A^-\Dq\iph).
\end{equation} 
Here $\Dq\imh = Q_i - Q_{i-1}$ is the jump across the interface at $i-1/2$,
for example.
For a linear system this is a generalization of the upwind method and is
first order accurate.

Second order accuracy is achieved by adding in correction fluxes:
\eql{wp1d2}
Q_i^{n+1} = Q_i^n - \DtDx (A^+\Dq + A^-\Dq) - \DtDx (\tF\iph - \tF\imh)
\end{equation} 
where
\eql{tF}
\tF\imh = \half\left( 1-\left|\frac{\lambda^p\Dt}{\Dx}\right|\right)
|\lambda^p|\W^p\imh
\end{equation} 
These terms convert the upwind method into a method of Lax-Wendroff type,
matching terms through $\Dt^2 A^2q_{xx}$ in the Taylor series expansion of the
solution at the end of the time step.  This method generates dispersive
errors, however, that can create large nonphysical oscillations near steep
gradients or discontinuities in a solution, such as shock waves.
To turn this into a ``high-resolution'' method, we use a wave limiter,
replacing $\W^p\imh$ in \eqn{tF} by $\tW^p\imh$, a limited version of the wave.
The wave $\W^p\imh$ is compared to the corresponding wave from the neighboring
Riemann problem, either $\W^p_{i-3/2}$ if $\lambda^p>0$ or $\W^p_{i+1/2}$ if
$\lambda^p<0$.  If the waves are of comparable magnitude the full correction
term is used for accuracy, but if there is a large discrepancy 
then the solution is not smooth at this point and a limited version is applied.
A wide variety of limiters have been developed, and we generally use the
MC limiter. 
See \cite{rjl:wpalg} or  Chapter 6 of  \cite{rjl_book} for
more complete details.

In two or three space dimensions the idea is the same, but now a 
1-dimensional Riemann problem must be solved normal to each edge or face of
the cell.  The resulting waves update the cell averages and correction
fluxes analogous to \eqn{tF} are used along with limiters in each direction.

In addition, to achieve second-order accuracy and good stability properties,
it is also necessary to use ``transverse Riemann solvers'' that further
modify the correction fluxes $\tF$ at each cell edge.
The method described above is based on propagating waves normal to each
interface.  In reality, the waves will propagate in a multidimensional
manner and affect cell averages in cells above and below those that are
directly adjacent to the interface.  

In two space dimensions, each ``fluctuation'' such as $A^-\Dq_{i-1/2,j}$ and $A^+\Dq_{i-1/2,j}$
that result from solving a Riemann problem in the $x$-direction
is split into two pieces using the eigenstructure of the
coefficient matrix $B$ in the $y$-direction, e.g.,
\eql{adqsplit}
A^+\Dq_{i-1/2,j} = B^-A^+\Dq_{i-1/2,j} + B^+A^+\Dq_{i-1/2,j}.
\end{equation} 
These two pieces will modify the correction flux at the edges $(i,j-1/2)$
and $(i,j+1/2)$ respectively to capture the transverse motion 
of the right-going wave.   Similarly, after solving a normal Riemann problem
in the $y$-direction using the $B$ matrix, transverse problems are solved
based on the eigenstructure of $A$.  The net effect of all these corrections
is to incorporate terms modeling the cross-derivative terms $BAq_{xy}$
and $ABq_{yx}$ of the Taylor series expansion in a properly upwinded manner.
More details can be found in 
\cite{rjl:wpalg} or  Chapter 21 of  \cite{rjl_book}.
The transverse correction terms are needed for accuracy, but
also have the effect
of improving the stability limit, allowing a Courant number near 1 to be used,
relative to the maximum wave speed in any direction.

In three space dimensions there are two transverse directions for each
normal Riemann solve, and terms modeling $CAq_{xz}$, etc. must also be
included.  Moreover, ``double transverse'' terms must be included, splitting
the result of a transverse solve into eigenvectors of the remaining
coefficeient matrix, and modeling terms such as $BCAq_{xzy}$.  The details
are presented in \cite{jol-rjl:3d} and fully implemented in Clawpack.

\subsection{The nonlinear fluid Riemann solver}\label{sec:nonlinRp}
The compressible fluid equations in Lagrangian form discussed in \Sec{lagrangian} 
can be reduced to the quasilinear form \eqn{quasilin} in which the Jacobian matrices
depend only on $q$ (for a spatially uniform fluid).   To apply the wave-propagation algorithm
we need to solve the Riemann problem orthogonal to each cell interface.  For nonlinear problems
this is usually done using an approximate Riemann solver, e.g., by replacing $f(q)_x$ by $\hat
Aq_x$, where the matrix $\hat A$ at each cell interface is chosen based on the data 
$Q_l$ and $Q_r$ to the left and right.
We use the f-wave formulation of the wave-propagation algorithm \cite{db-rjl-sm-jr:vcflux}, in
which the jump in flux $ f(Q_r) - f(Q_l)$ is split into eigenvectors of an approximate Jacobian
matrix, rather than the jump in $Q$.  This leads to an algorithm that is conservative for any
choice of approximate Jacobian and also extends naturally to the case of spatially varying
fluxes, as required near the fluid-solid boundary and discussed further below.

Rather than choose an approximate Jacobian $\hat A$ and then determining its eigenvectors and
eigenvalues, we simply choose the set of eigenvectors and associated wave speeds based on the
data and wave forms expected to result from this data.  These vectors form a matrix $\hat R$ and
we then solve
$\hat R \beta = f(Q_r) - f(Q_l)$ for the vector of wave strengths $\beta$.
The choice of vectors in $\hat R$ and associated wave speeds $\hat \lambda$ implicitly defines
the Jacobian approximation $\hat A = \hat R\hat\Lambda\hat R\inv$, but this matrix is never
needed.

The eigenvectors are taken to be the vectors displayed in \eqn{r1-6} and \eqn{r789}.  Recall that in
the fluid case there are only two nonzero eigenvalues corresponding to the first two
eigenvectors.  For the eigenvector corresponding to  $\lambda^1 <0$ we use
$\lambda^1 = -\sigma^{11}_{\epsilon^{11}}$ evaluated in the left state $Q_l$, while the eigenvector
corresponding to $\lambda^2 >0$ is determined using
$\lambda^2 = \sigma^{11}_{\epsilon^{11}}$ evaluated in the right  state $Q_r$.
These vectors have nonzero components only in positions 1 and 7 and so the values
of $\beta^1$ and $\beta^2$ can be determined by solving a $2\times 2$ system:
\eql{beta12}
\bcm 1&1 \\ 
\rho_l\lambda_l^1 & \rho_r \lambda_r^2 \ecm
\bcm \beta^1\\ \beta^2\ecm  
= \bcm \Delta f^1 \\ \Delta f^7 \ecm.
\end{equation} 
The solution is
\eql{betas2}
\beta^1 = \frac{\rho_r \lambda_r^2\Delta f^1  - \Delta f^7}{\rho_r \lambda_r^2
- \rho_l \lambda_l^1}, \qquad \beta^2=\frac{\Delta f^7 - \rho_l \lambda_l^1 \Delta f^1}{\rho_r \lambda_r^2 - \rho_l \lambda_l^1}.
\end{equation} 
The remaining waves do not propagate and do not come into the wave-propagation algorithm.

\subsection{The linear elastic Riemann solver}\label{sec:hetero}

In the linear elastic material modeling bone, we take a similar approach and again use the f-wave
formulation of the algorithm.  In this case there
are six waves with nonzero wave speeds given by the eigenvectors in \eqn{r1-6}.  The
eigenvectors are independent of $q$ in the linear case, but can be spatially varying
to represent varying bone structure, so the coefficients $C_{ij}$ in \eqn{hookes} can vary from
one grid cell to the next.  Similar to the nonlinear case described above, to compute the
decomposition of the flux difference into propagating waves we define the three left-going 
eigenvectors $r^{1,3,5}$ (with the minus sign in \eqn{r1-6}) based on the coefficients in the
left state, while the right-going eigenvectors $r^{2,4,6}$ are defined using the coefficients in
the right state.  Note that the flux vector $f(q)$ from \eqn{qfgh}, and hence any jump in flux,
has zeros in three components which are easily seen to lead to $\beta^7=\beta^8=\beta^9$ when the flux
difference is written as a linear combination of the eigenvectors, and the six remaining components of the
flux difference uniquely define the coefficients $\beta^1$ through $\beta^6$ for the six propagating
waves.  $\beta^1$ and $\beta^2$ are the same as \eqn{betas2}, the other weights are
\eql{betas34}
\beta^3 = \frac{\rho_r \lambda_r^4\Delta f^4  - \Delta f^8}{\rho_r \lambda_r^4
- \rho_l \lambda_l^3}, \qquad \beta^4=\frac{\Delta f^8 - \rho_l \lambda_l^3 \Delta f^4}{\rho_r \lambda_r^4 - \rho_l \lambda_l^4},
\end{equation}
\begin{equation*}
\beta^5 = \frac{\rho_r \lambda_r^6\Delta f^6  - \Delta f^9}{\rho_r \lambda_r^6
- \rho_l \lambda_l^5}, \qquad \beta^6=\frac{\Delta f^9 - \rho_l \lambda_l^5 \Delta f^6}{\rho_r \lambda_r^6 - \rho_l \lambda_l^5}.
\end{equation*} 

\subsection{Interfaces and the Cartesian Grid}
\label{sec:interfaces}
\begin{figure}
\begin{center}
a)\includegraphics[height=1.8in]{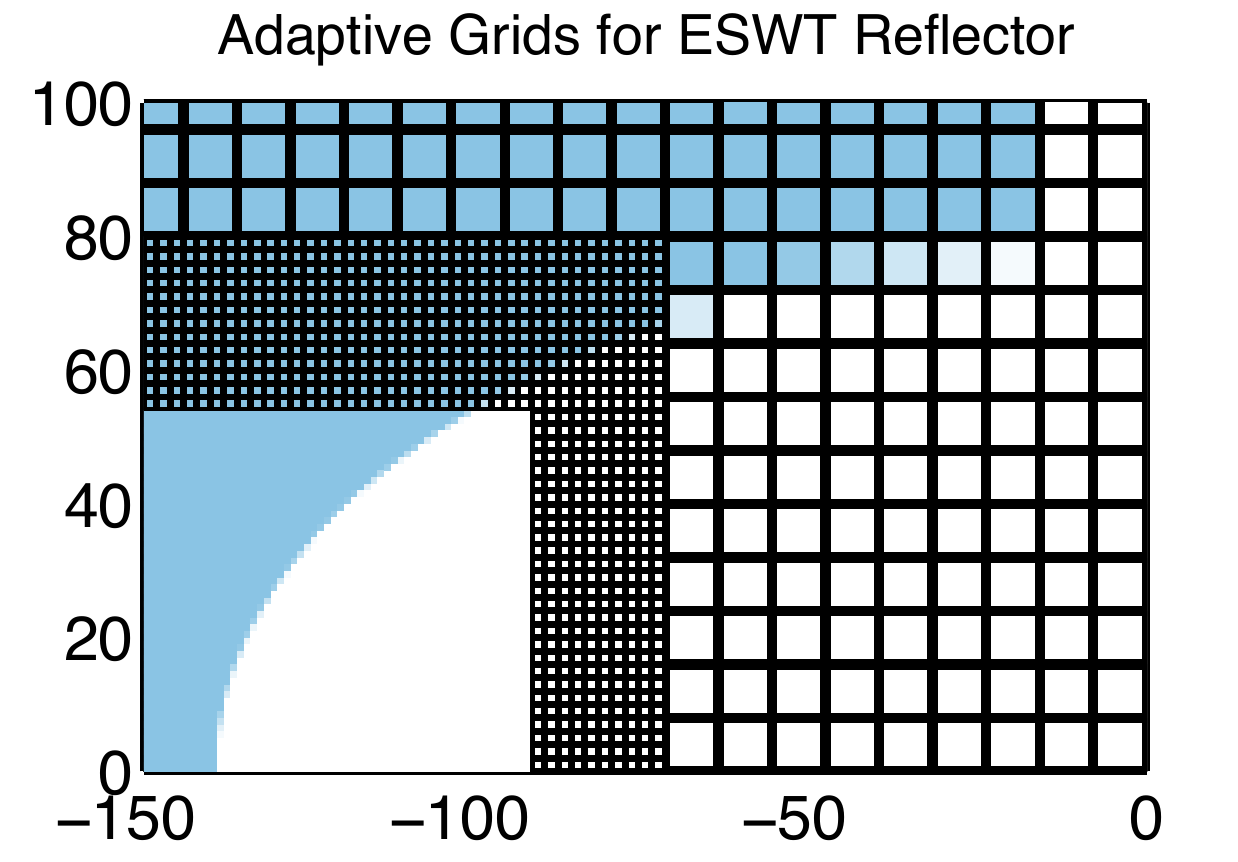}
b)\includegraphics[height=1.75in]{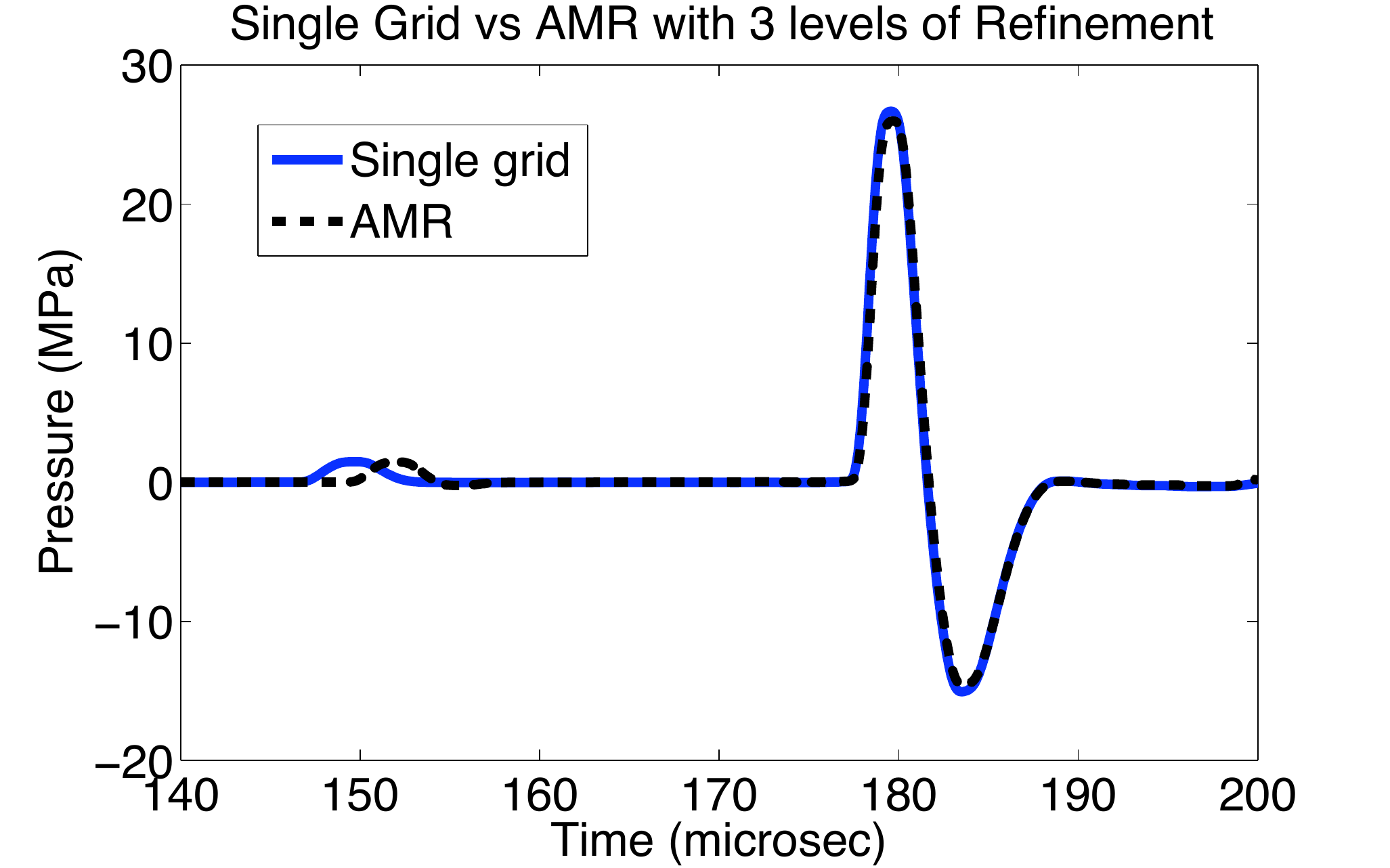}
\caption{ a) Resolution of the ellipsoid reflector with different levels of AMR.  b) 2d axisymmetric calculation.  The second is a comparison of the waveform obtained using AMR and a uniform grid.  The finest grid resolution in the AMR calculation is the same as the resolution on the uniform grid. }
\label{fig:amr_vs_singlegrid}
\end{center}
\end{figure}

In ESWT the pressure wave must propagate through a variety of materials, and in general the interfaces between different materials do not align with the grid.  In our calculations we use a Cartesian grid.  To handle grid cells that contain two materials, we perform a weighted average of the material properties.  The stress-strain relationship in the averaged grid cells is taken to be that from linear elasticity, even if one of the materials is fluid.  This approach is feasible because we use AMR to refine around the interfaces between the two materials.  By using a fine enough grid, we are able to reduce the error introduced by the weighted average approximation.  Figure \ref{fig:amr_vs_singlegrid} a) illustrates the interface between the fluid and the brass reflector from an axisymmetric calculation.  Three grid resolutions are shown in this figure: a coarse grid on the right, a level 2 grid that is refined by a factor of 4 in each direction in the middle, and the finest grid on the left, where the grid lines are not drawn.

Figure \ref{fig:amr_vs_singlegrid} b) shows a comparison between the pressure wave measured at F2 for an AMR calculation versus a single grid calculation.  The single grid calculation took 269 minutes to complete, just over 6 times as long as long as the run using AMR which finished in 44 minutes.  These calculations were performed serially on a 2.8 GHz dual core AMD Opteron machine with 32 GB of memory.  It's clear that the two calculations yield comparable pressure waves.  The biggest difference is is in the direct wave arriving around $t=150$, which is not being resolved in the AMR calculation because we have refined only in the vicinity of the reflected wave of primary interest. 

\subsection{Adaptive Mesh Refinement}
\label{sec:amr}
\begin{figure}[h!]
\begin{center}
\includegraphics[height=2.5in]{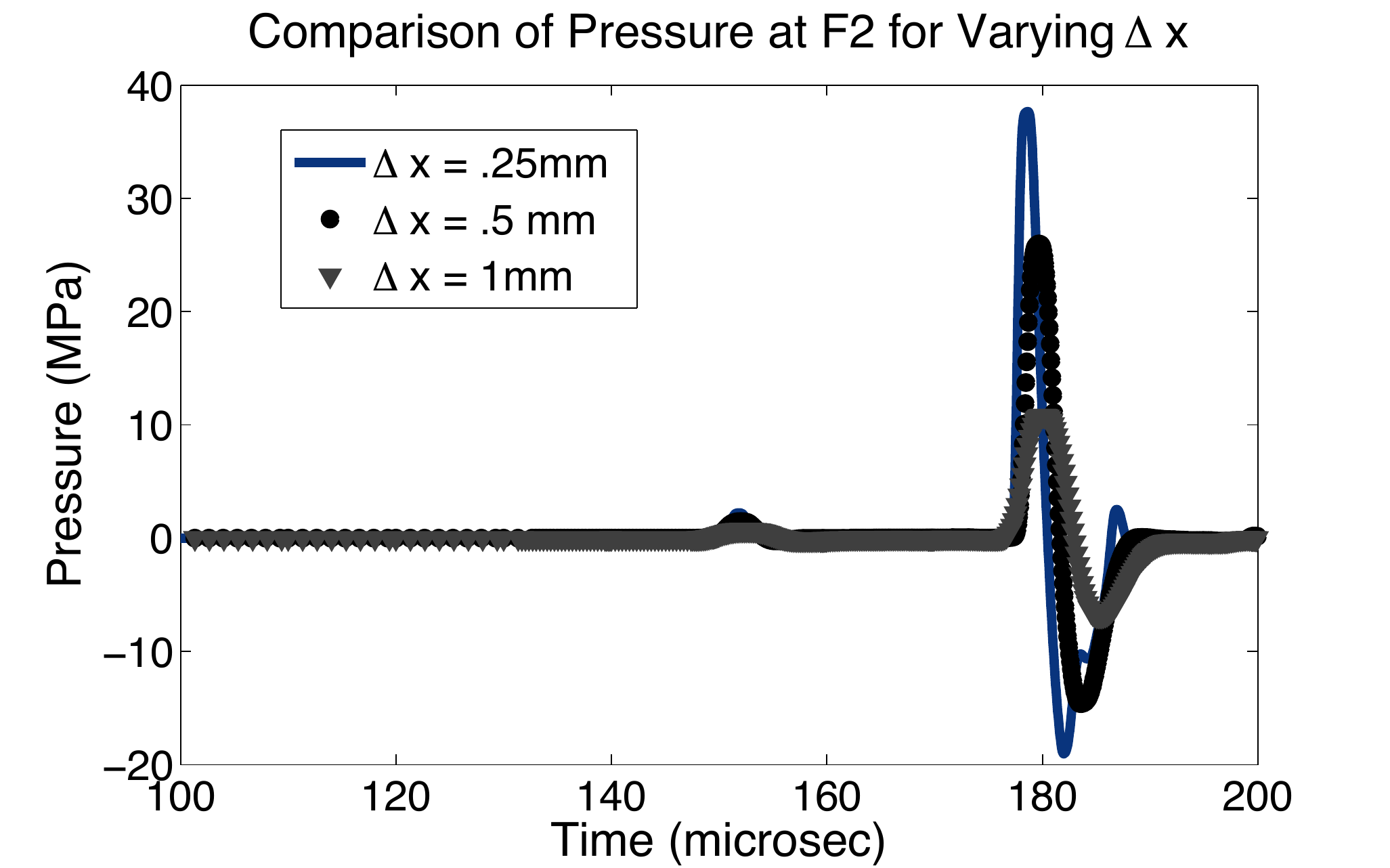}
\label{fig:steepening}
\caption{Effect of grid size on shock wave profile.  As the grid is refined for the same initial condition, the shock wave profile steepens. The solution eventually converges to a profile with the same magnitude, though the convergence rate is only first order near a discontinuity. }
\end{center}
\end{figure}

The pressure waveform found in ESWT contains a very thin region of high pressure that can not be 
resolved without a highly refined mesh.  In Figure \ref{fig:steepening} we investigated the effect of grid 
refinement on the shock wave profile and found that with grid resolution greater than $0.25$mm, the 
wave form at F2 was not a shock.  Note that near the shock we only expect our method to be first order, but the solution does converge to a shock as the grid is refined.  Our calculations are done with the adaptive mesh refinement (AMR) in the style of Collela-Berger-Oliger \cite{mjb-col:amr}, \cite{mjb-ol:amr}.  
The AMR algorithms used in Clawpack are more fully described in \cite{mjb-rjl:amrclaw}.  For the three-dimensional calculations, a similar AMR algorithm is used, as implemented in Chombo.  Here we only briefly review the main ideas.

The computational domain is covered by a rectangular level 1 grid, typically at a coarse resolution.  Rectangular patches of the grid may be  covered by level 2 grids, refined by some specified refinement ratio in each direction.  
Since we use explicit methods, the CFL condition generally requires that the time step be refined by the same factor on the level 2 grids, so several time steps must be taken on each level 2 grid for each time step on the level 1 grid.  The level 1 grid is advanced first, and for each time step on the level 2 grid, ghost cell values around the boundary are filled in either by copying from adjacent grids at the same level, or using space-time interpolation from the level 1 grid for ghost cells that do not lie in an adjacent grid.  This entire procedure is repeated recursively to obtain higher levels of refinement, e.g. some portion of the collection of level 2 grids may be covered by level 3 grids and so on.

In order to adaptively refine the grid, it is important to specify appropriate refinement criteria.  The perturbations to the strain are small, so gradients in the strain are too small to use as reliable refinement criteria.  However, the small strains result in large changes in the pressure, so we refine in the area near the pressure wave.  In order to handle the interfaces between two materials, we also use large gradients in background density as a secondary refinement criterion.  
Cells that are flagged as needing refinement are clustered into rectangular patches using the algorithm of Berger and Rigoutsis \cite{mjb-rig:cluster}.  Regridding is done every few steps on each grid level in order to track propagating waves.  
Regions are automatically de-refined once the wave passes by, since in these regions cells are no longer flagged as needing refinement.

Figure \ref{fig:steepening} illustrates the behavior of the ESWT waveform as the grid is refined.  What is evident from these experiments is that a coarse grid will not effectively capture the development of the shock, so around the propagating wave, we need at least $\Delta x = .25mm$ resolution.  As the wave steepens into a shock, we no longer expect second order convergence, because in the region around a discontinuity, our methods are first order.  However, since the discontinuities occur in a small region of the domain, the overall methodology is still second order.  

In order to efficiently calculate a reasonable ESWT waveform in three dimensions, we utilized ChomboClaw 
\cite{chomboclaw}, which uses the adaptive mesh refinement routines of CHOMBO with the wave 
propagation solvers of Clawpack.  This code can be run in parallel using MPI on an NSF
TeraGrid computer at TACC.  
We found that the ChomboClaw code scaled up to 128 processors: the same experiment ran approximately 4 times faster on 128 processors than on 32.  
However, we need to do further testing on larger systems in order to validate that ChomboClaw will scale as well as CHOMBO,
which is reported to scale well to  10,000 processors
\cite{colella_communication}.  

\section{Results}
\label{sec:results}
We have used the approach described above to model ESWT pressure waves interacting with three-dimensional bone 
geometries comprised of idealized materials.  We have modeled both simple objects 
that have been used in laboratory experiments \cite{amath_apl_nonunion}
as well as complex three-dimensional geometries
extracted from CT scans of patient data \cite{kfagnan_mchang_ho}.  Here we present 
results that demonstrate the efficacy of the Lagrangian formulation, as well as examples of calculations 
performed using real three-dimensional geometries.  

The calculations were initialized using pressure data obtained from a 2D axisymmetric calculation where we modeled the full geometry of the ellipsoidal reflector.  The reflector was modeled using linear elasticity with material properties that can be found in \cite{fagnan:phd}.  We assumed the fluid was water with the corresponding parameters for the Tait equation of state found in \Sec{euler}.   We saved the data at $t=116\mu s$ and used this to restart future calculations.  For the three-dimensional initial condition, we rotated the 2D data about the $x$-axis.  The material properties of averaged bone were obtained from \cite{martin_burr_sharkey} and used in the heterotopic ossification, cylinder and sphere.  

We have found in our experiments that interfaces between materials with large impedance differences have the most significant effect on maximum stress and energy deposition.

\subsection{Reflection and focusing}
\label{sec:focus}
\begin{figure}[h!]
\begin{center}
\includegraphics[width=.45\textwidth]{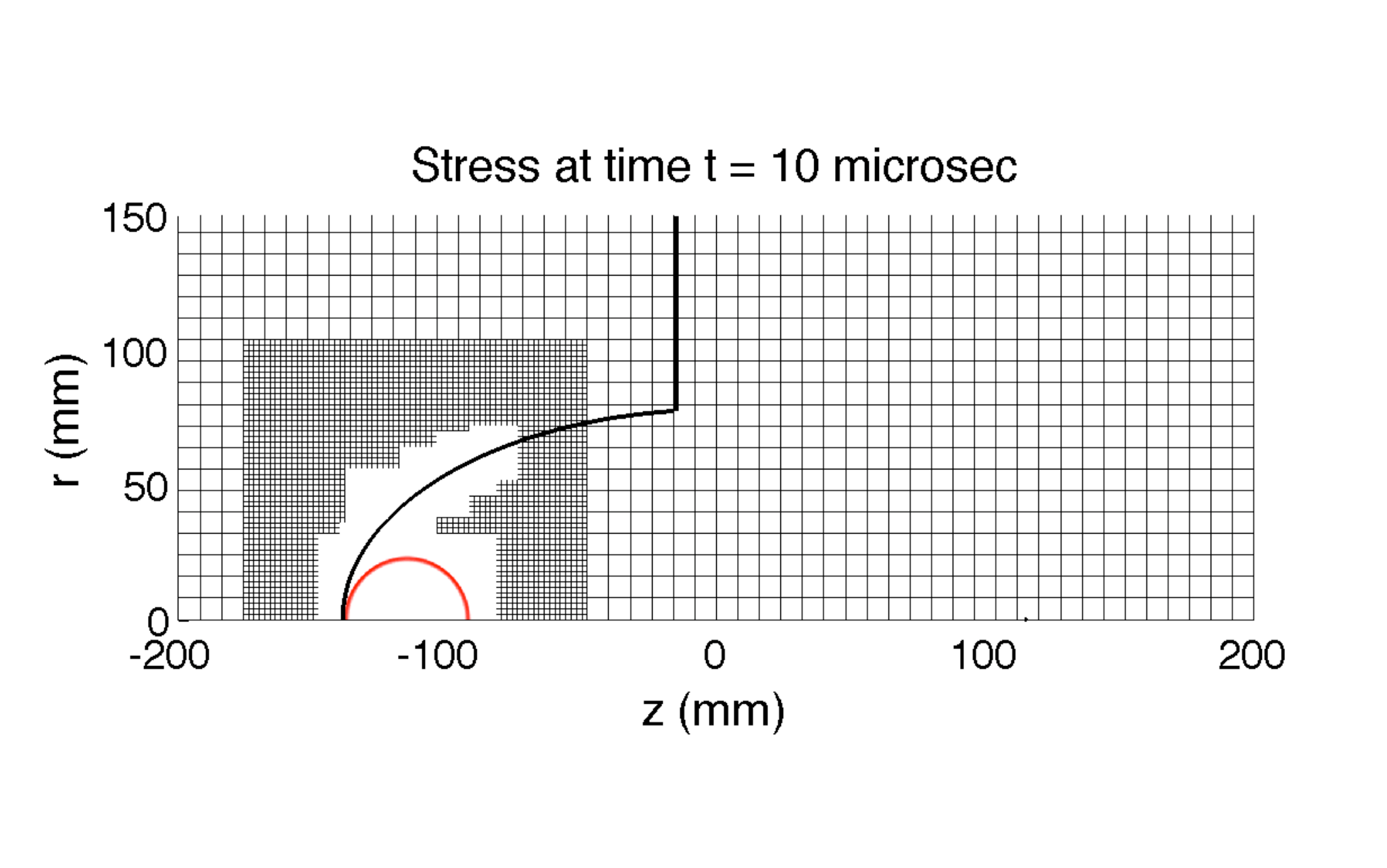}
\includegraphics[width=.45\textwidth]{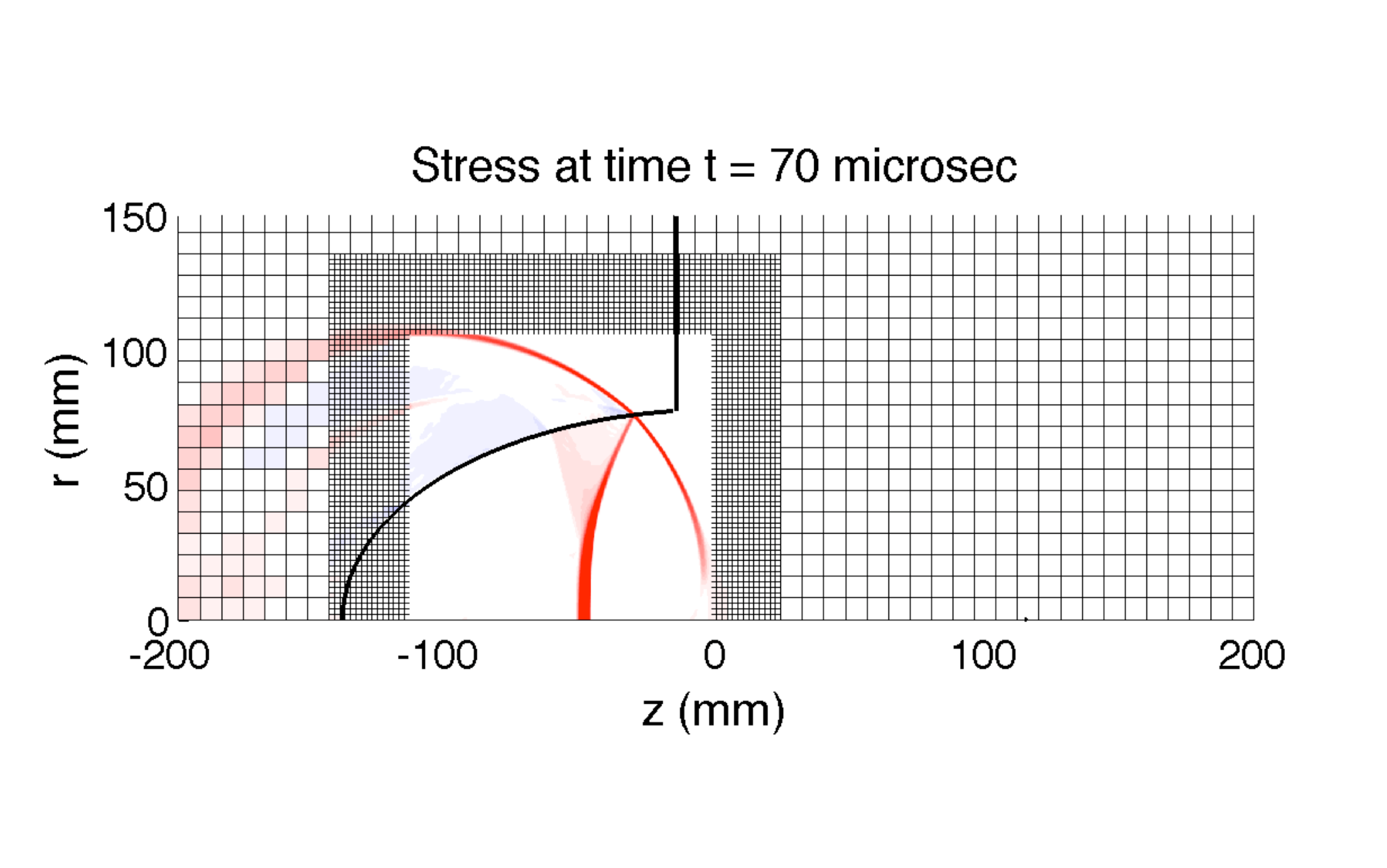}\vspace{1mm}
\includegraphics[width=.45\textwidth]{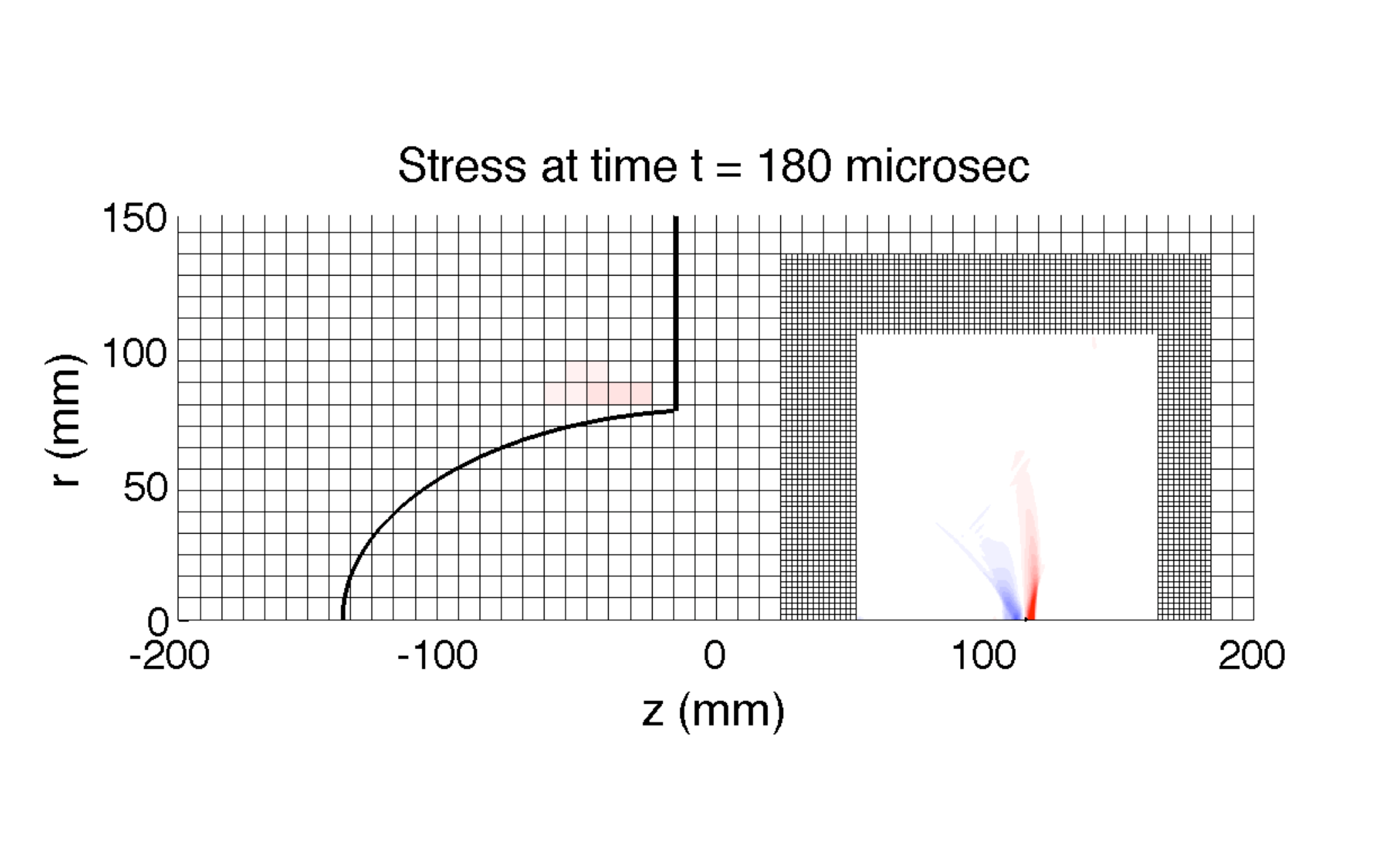}
\includegraphics[width=.45\textwidth]{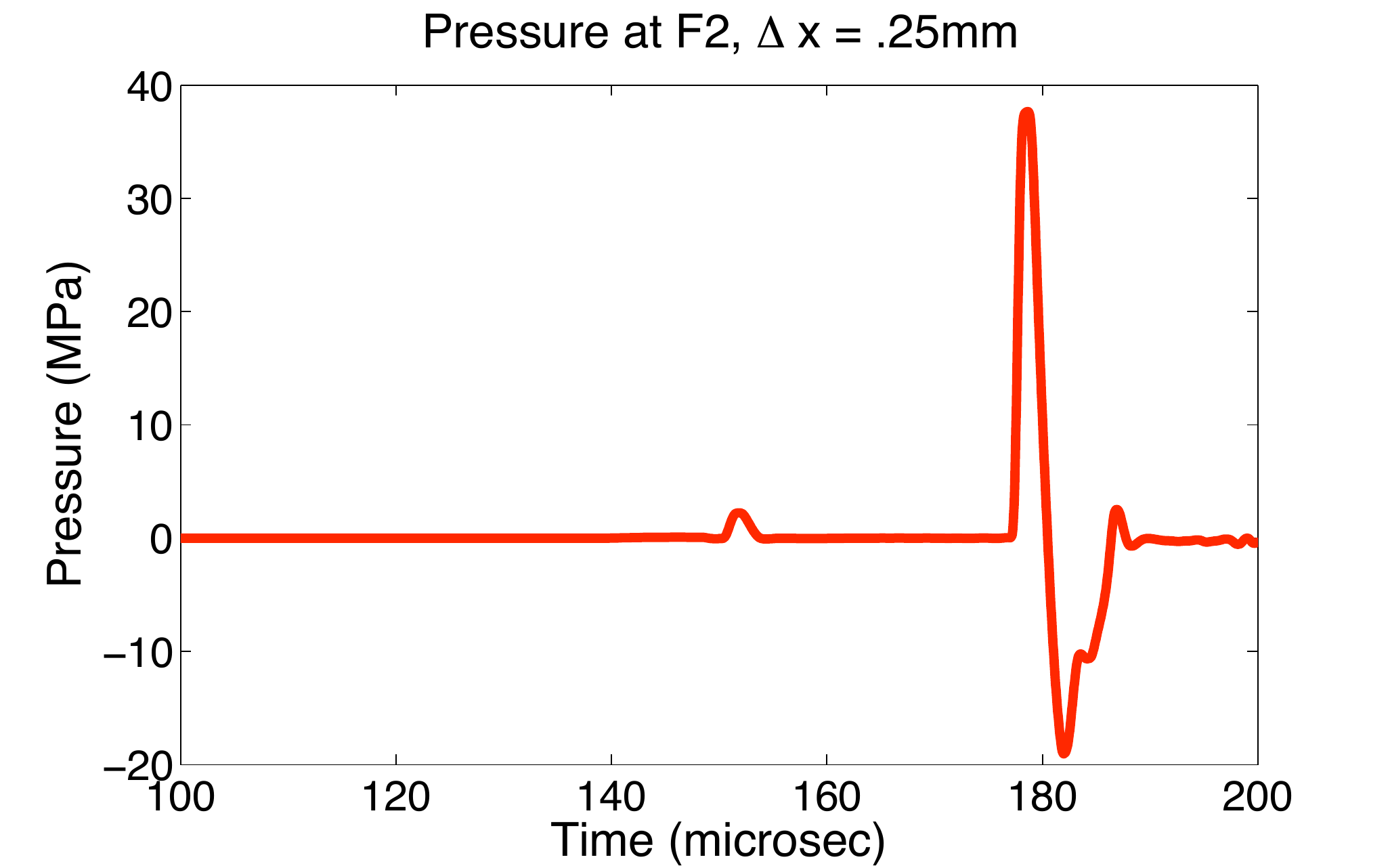}
\caption{
\label{fig:amr_interface}
Axisymmetric calculation of the pressure pulse generated by a spherical high-pressure bubble
centered $z=-115$ (the focus F1 of the ellipsoidal reflector). Three levels of AMR are used and
grid lines are shown only on levels 1 and 2. The level 3 grid has a resolution of $\Delta z =
\Delta r = 0.25$mm. a) At $t=10$ the pulse has nearly reached the reflector. b) At $t=70$ the
incident, transmitted, and reflected pulses are visible. c) At $t=180$ the reflected pulse has
focussed near $z=115$ (the focus F2). d) The time history of the pressure at F2. The direct
(unreflected) wave passes F2 at $t\approx 150$ and the focussed pulse arrives at $t\approx 180$.}
\end{center}
\end{figure}

In Figure \ref{fig:amr_interface}, we show an axisymmetric calculation of the ESWT wave propagation
and focusing in water alone, in a domain bounded by the ellipsoidal reflector of the Dornier HM3.
Figure \ref{fig:amr_interface} a) shows the initial spherical propagation of the pressure wave, as
well as the grids where the calculation is being refined. The grid must be refined around the
pressure wave as well as the reflector. Figures \ref{fig:amr_interface} b) and c) show the
propagation of the wave and evolution of the the adaptive grid structures. At later times the grid
is only being refined near the pressure wave. The sharp results and absence of spurious
oscillations in the pressure measurement at F2 indicate that AMR together with our Cartesian grid
approach enables us to capture the reflection at the interface.
\clearpage

\subsection{Axisymmetric sphere}
\label{sec:sphere}
We used an axisymmetric test problem in order to compare the solutions obtained with the two-dimensional and three-dimensional codes.  The initial condition for this experiment was an analytic form for an ESWT pressure wave used in \cite{bailey_oleg}.  In the two-dimensional case, we specified the pressure as a function of the radial distance from F1(-115,0).  In the three-dimensional case, we rotated the same two-dimensional initial condition about the $x$-axis.  The grid resolution was $\Delta x = .25 mm$.  

\begin{figure}[h!]
\begin{center}
\includegraphics[width=2.5in]{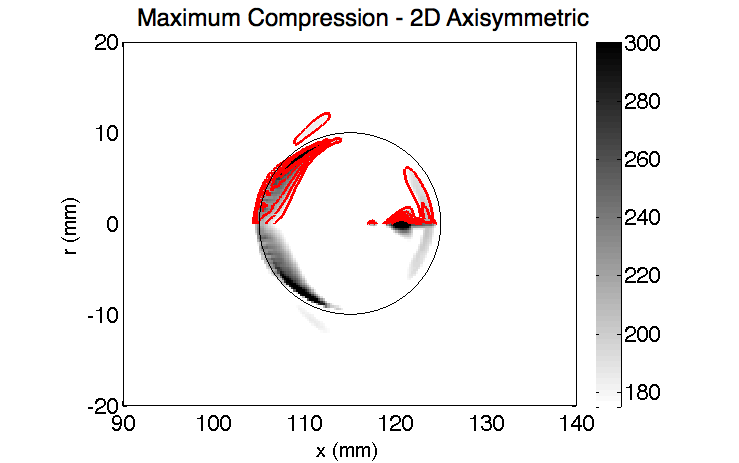}\hspace{-10mm}
\includegraphics[width=2.25in]{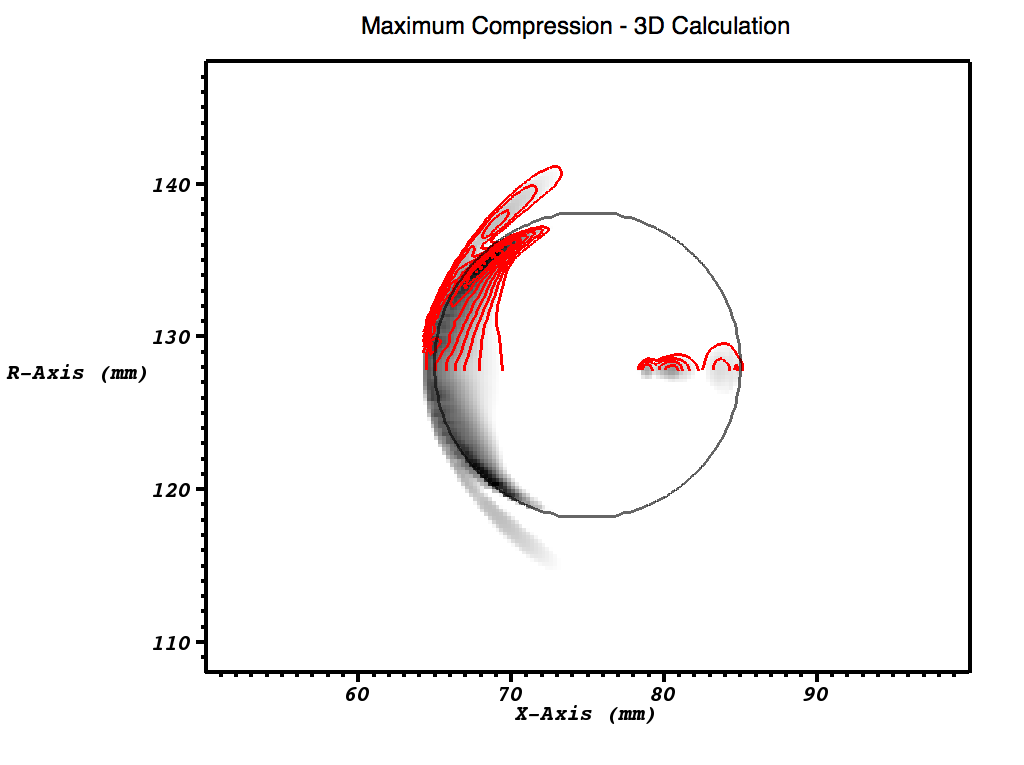}
\vskip 5pt

\includegraphics[width=2.5in]{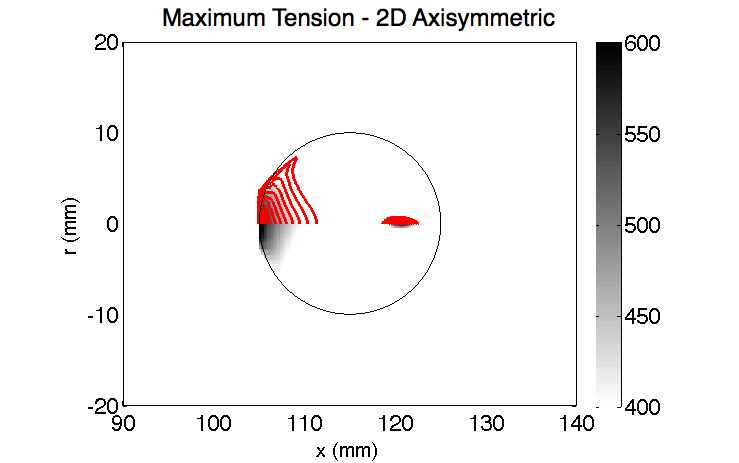}\hspace{-10mm}
\includegraphics[width=2.25in]{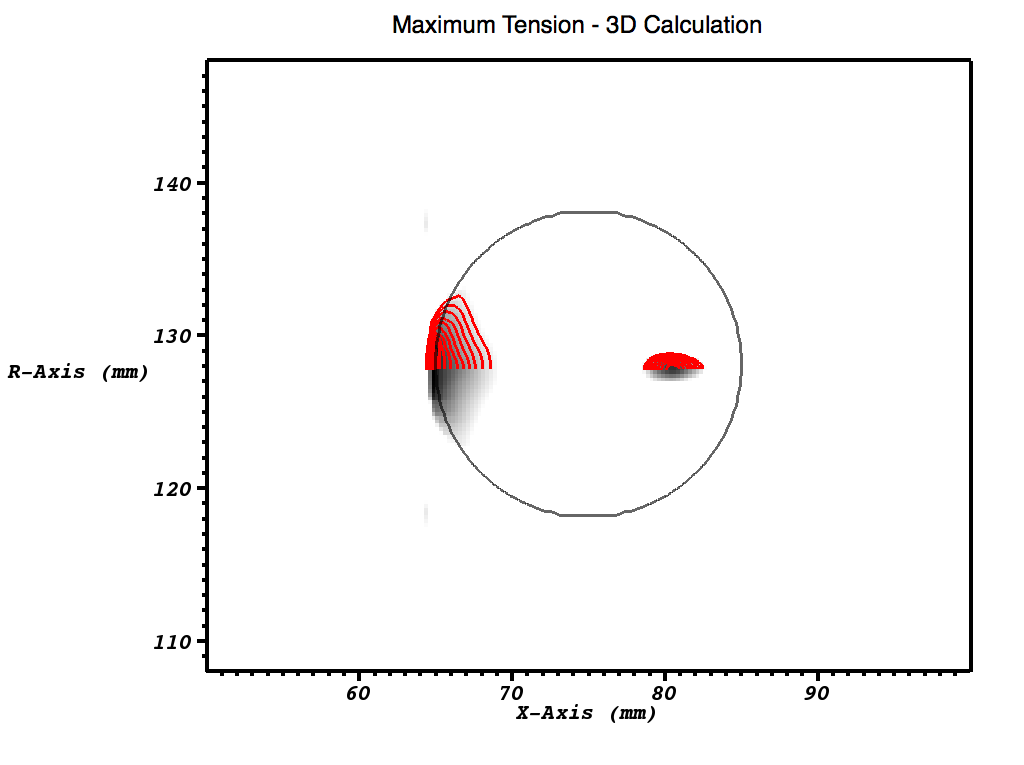}

\vskip 5pt
\includegraphics[width=2.5in]{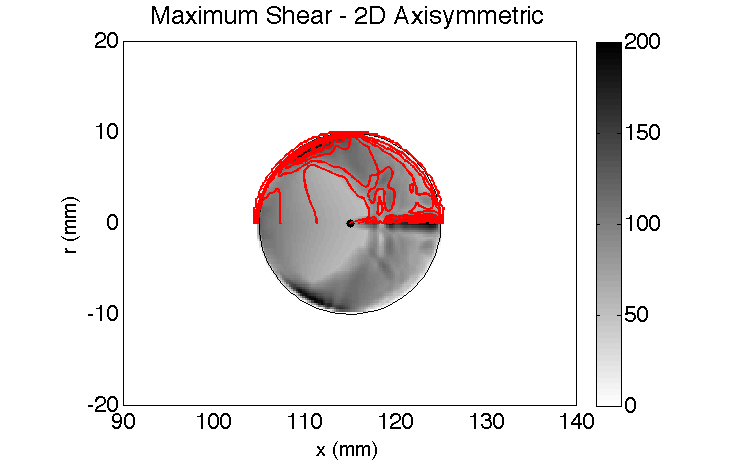}
\includegraphics[width=2.25in]{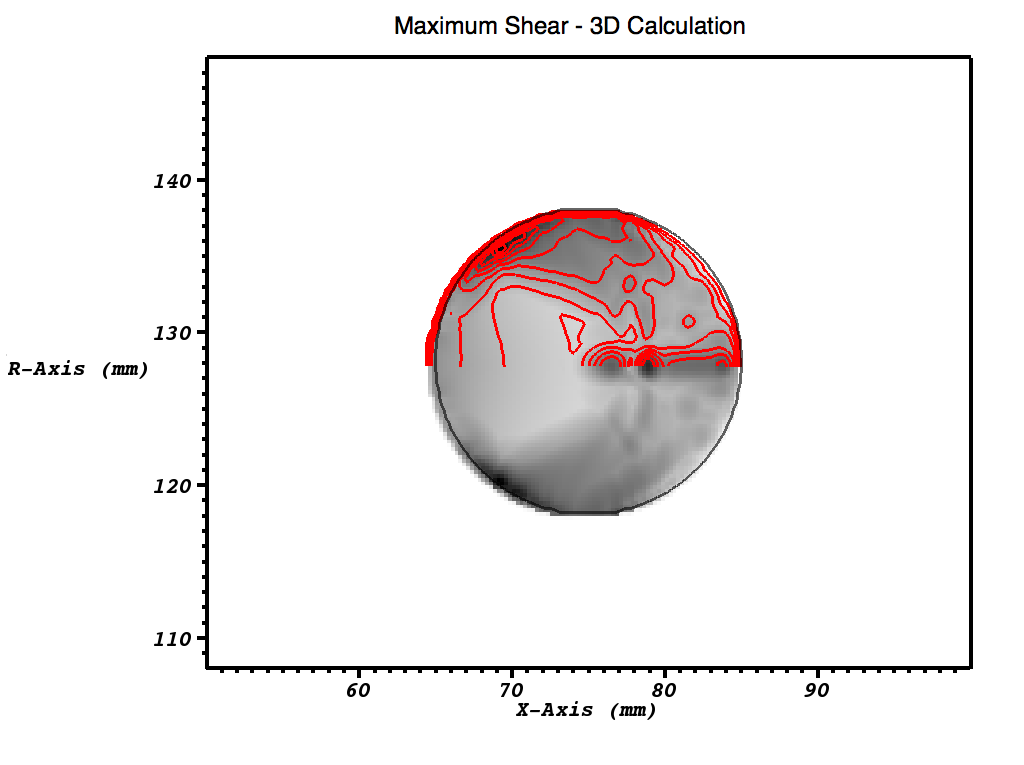}
\caption{Results from 
calculation of a shockwave interacting with an acrylic sphere.
The left column shows 2D axisymmetric results and the right column shows
a corresponding cross section of full 3D calculation. 
Top: Maximum compression.  Middle: Maximum tension.  Bottom: Maximum shear.
\label{fig:sphere_comp_2d_3d}}
\end{center}
\end{figure}

\begin{figure}[h!]
\begin{center}
\includegraphics[width=5.in]{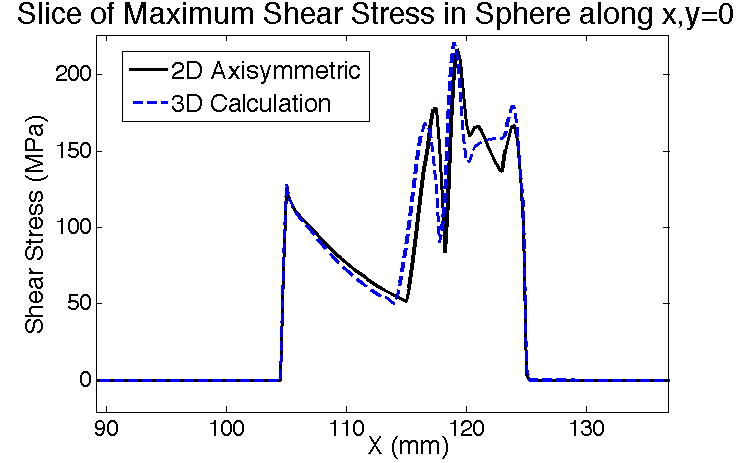}
\caption{Comparison of maximum shear stress from 2D and 3D calculations
as a function of $x$ along $y=z=0$.  The difference in the results is likely caused by averaging of the initial 
condition onto the 3D domain and the boundary conditions on the axisymmetric calculation at $r=0$, but 
the two calculations predict comparable location and magnitude of maximal shear stress deposition. 
\label{fig:sphere_comp_shear}}
\end{center}
\end{figure}

Results with contour lines are in Figure \ref{fig:sphere_comp_2d_3d}.  The
maximum values in each of the three cases are nearly the same, but there are
slight discrepancies in the contour lines.  Figure
\ref{fig:sphere_comp_shear} shows a 1D slice of the maximum shear
calculation in the 2D and 3D codes, which makes it clear that the peak of maximal shear stress is in the same location and has the same value.  The general shape of the maximum stress deposition pattern are similar in both cases.  The difference in the two solutions is likely caused by the solid wall boundary condition that is used at $r=0$.  Only waves that are propagating normal to that boundary are perfectly reflected, otherwise some error is generated.  


\subsection{Nonunions}
\label{sec:cylinders}
ESWT has recently been used in the treatment of non-unions or bone
fractures that fail to heal \cite{eswt_bone}.  One question that
is of interest to clinicians is whether or not the angle of treatment
has an effect on healing.  We assume that healing is related to the
magnitude of stress  applied near the treatment area, although the
connection between the applied force and biological response is not
yet understood.  
In the fluid there is no shear stress.
However, at the liquid-solid interfaces, shear stresses are generated
by the shockwave and stimulate motion both at the surface and within
the material.  The motion of the biological materials (e.g. the
periosteum, interstitial fluid, mechanotransduction) is likely to be important in the
healing process \cite{frangos, prendergast1997, weinbaum1994, claes1999, isaksson2006,park1998,turner1998}, and modeling
the magnitude and location of the stress deposition is a good first
step toward understanding the shear and tensile displacements caused by ESWT.  
We should stress, however, that the healing mechanisms are not well
understood and we are not claiming that magnitude of the applied
stress is the most important or only biological mechanism involved
in the healing process.  As mentioned in \Sec{intro}, several studies
have indicated that cyclic application of mechanical loading leads to the generation
of new bone.  The work of Isaksson, et. al. \cite{isaksson2006}, indicates that the most accurate
predictors for bone healing are those based on shear strain and fluid flow, however, there is no 
single model that can predict all features of the healing process, so more work is necessary \cite{morgan2008}.

\begin{figure}[h!]
\includegraphics[height=1.5in]{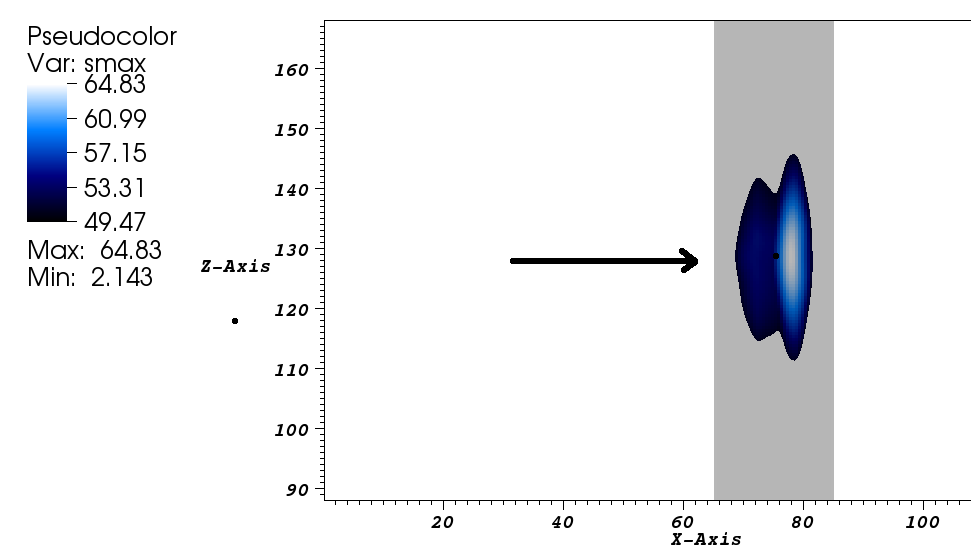}
\includegraphics[height=1.5in]{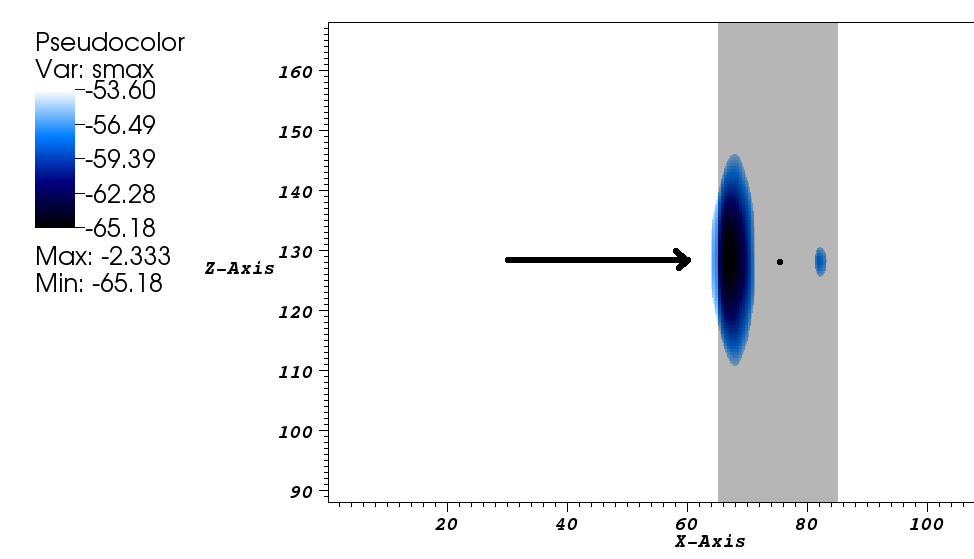}
\includegraphics[height=1.5in]{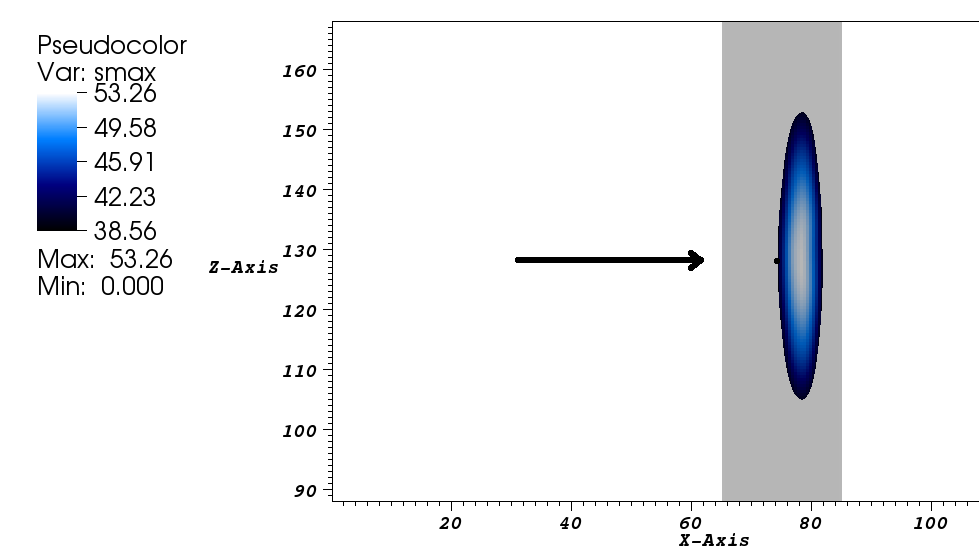}
\caption{Three-dimensional results for the direct treatment of a complete cylinder.  This figure shows 2D slices of maximum compression, tension and shear along y=0 for treatment where the ESWT wave propagates along the x-axis, as indicated by the arrow.  The dot illustrates the location of F2.}
\label{fig:90cylin}
\end{figure}

\begin{figure}[h!]
\includegraphics[height=1.5in]{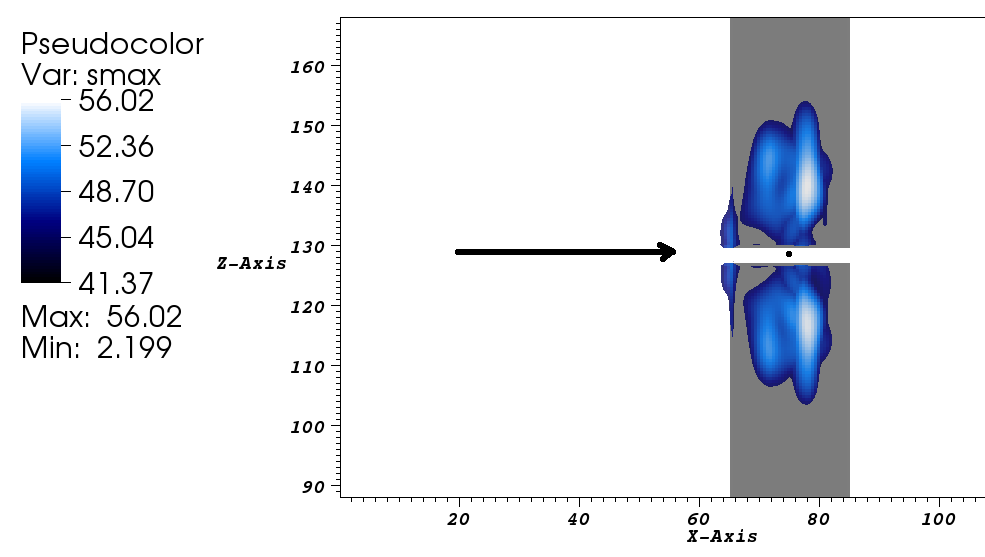}
\includegraphics[height=1.5in]{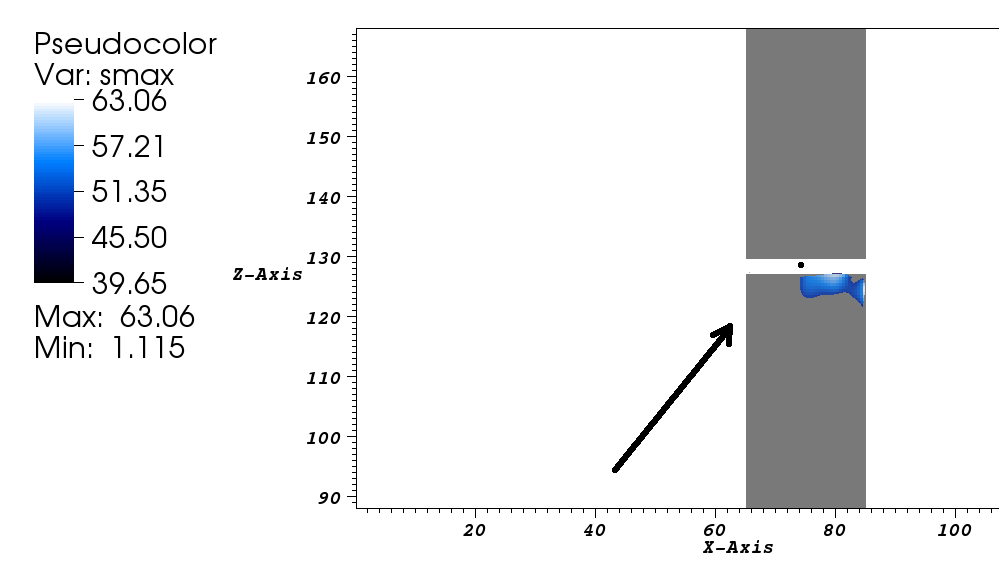} 
\includegraphics[height=1.5in]{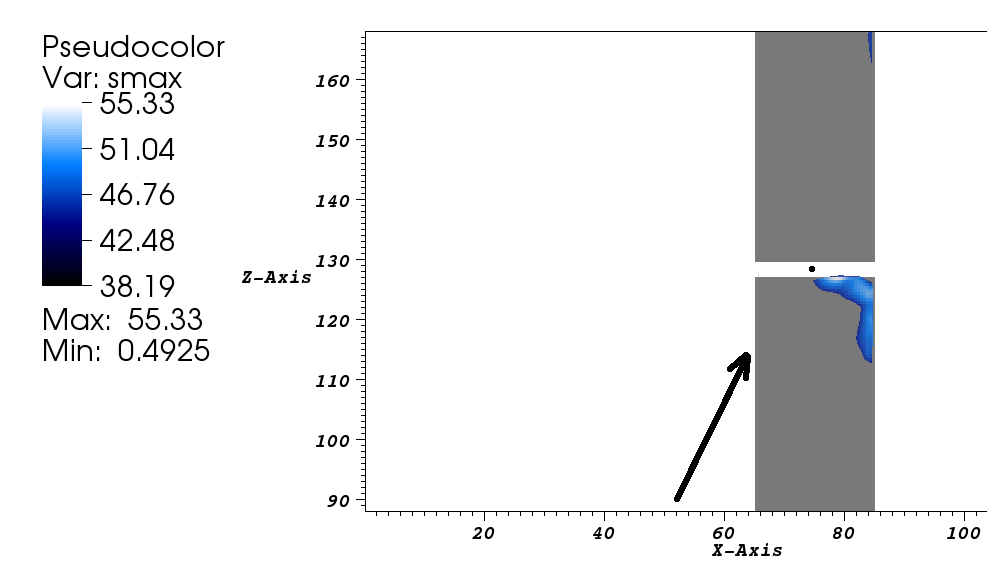}
\caption{Three-dimensional results for the direct treatment of a broken cylinder.  This figure shows 2D slices of maximum compression along y=0 for treatment along the x-axis, 45 degree rotation and 60 degree rotation about the y-axis.  The arrows indicate the angle of treatment in each case and the dot illustrates the location of F2.}
\label{fig:90gapcylin}
\end{figure}

In an actual treatment, the clinician generally sets up the device so that the focus is aligned with the ailment.  For example, in the case of a broken bone, the clinician will set up the device so that F2 is in the center of the break.  However, given the heterogeneous media, it is not clear that the maximal stresses will actually be observed at F2, as would be expected in pure water.  We used our model to investigate the location of maximal stress deposition relative to F2.  In these calculations we considered two different geometries, a complete cylinder, representing the long shaft of a healthy bone, and a broken cylinder, representing a nonunion.  The results from calculations where the idealized bone was perpendicular to the direction of the pressure wave front are shown in figures \ref{fig:90cylin} and \ref{fig:90gapcylin}.  We found that the break has a significant impact on the location of stress deposition.  

We used these geometries to perform a variety of experiments.   We rotated the direction of treatment by 45 and 60 degrees relative to the $x$-axis and calculated both the magnitude of the maximum compressive, tensile and shear stresses, as well as the distance from the focus F2 of the device. \ignore{Is x-axis defined?}

In the case of the broken cylinder, the maximum stress deposition in the direct experiment is similar to that of the unbroken cylinder, except that the there are two locations of maximal stress deposition on either side of the break.  The pressures in the bone are larger than in the fluid due to reflection at the fluid-solid interface, so the contours of maximum stress are concentrated on either side of the gap.  The location along the x-axis is nearly the same as in the unbroken cylinder, and the distances from the ideal focal point, F2, are also similar. 

As the angle of treatment is varied, there is less of a shift in the z-direction for the shear and compressive stresses.  This is caused by the impedance difference between the fluid and solid material at the gap, which is located close to F2.  If the gap were shifted along the z-axis from the focal point, there would be a corresponding shift in the location of maximum shear and compression.  Geometrically, the shape of the regions of compressive and shear stress are quite different from the direct case.  Instead of being an ellipsoidal shape, the regions are compressed into the corner of the lower-half of the cylinder.  Again, this is caused by the impedance jump at the fluid-solid interface.  The region of maximum tension deposition is similar to that of the unbroken cylinder case, though it is also affected by the gap and the tension is concentrated on the upper half of the cylinder.

It is clear from the literature
\cite{park1998,prendergast1997,lacroix2002,claes1999,isaksson2006,carter1998,goodship1985},
that mechanical loading is important in bone healing. The implication of our
computational experiments is that the angle of treatment will affect stress
deposition and therefore may be important in treatment optimization.  For
example, in order to maximize shear stress at the tissue-bone interface, our
preliminary computations indicate that it might be best to treat the patient at an oblique angle.  However, if the goal is to maximize shear stress in the gap of the broken bone, then treating the patient at a 90 degree angle may be better than treating at either the 45 or 60 degree angle.  We stress however that the biological mechanisms must be better understood and more experiments must be done in conjunction with laboratory and clinical treatments before these calculations could be used to make specific clinical recommendations.

\begin{figure}[h!]
\begin{center}
\includegraphics[height=2.5in]{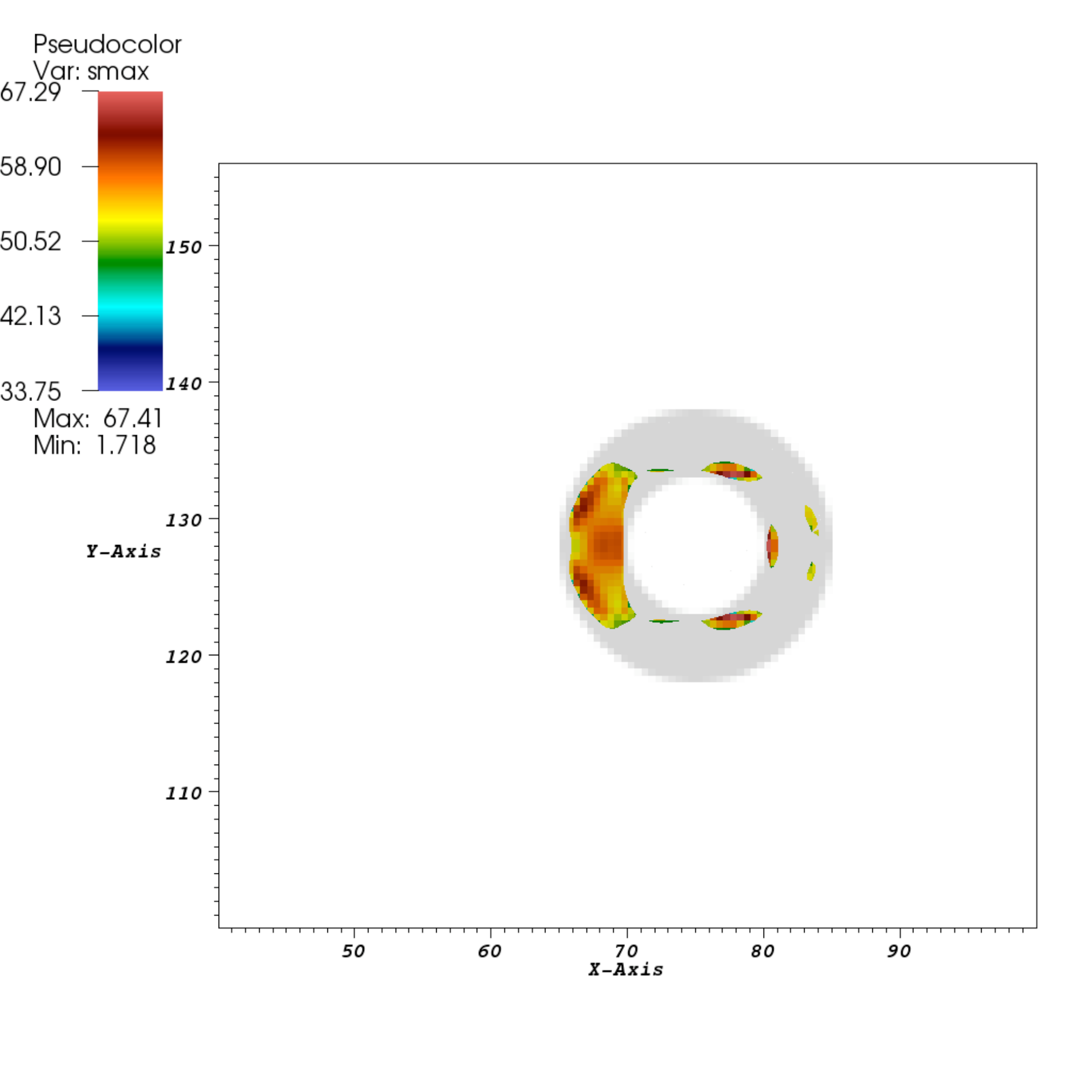}\hspace{5mm}
\includegraphics[height=2.5in]{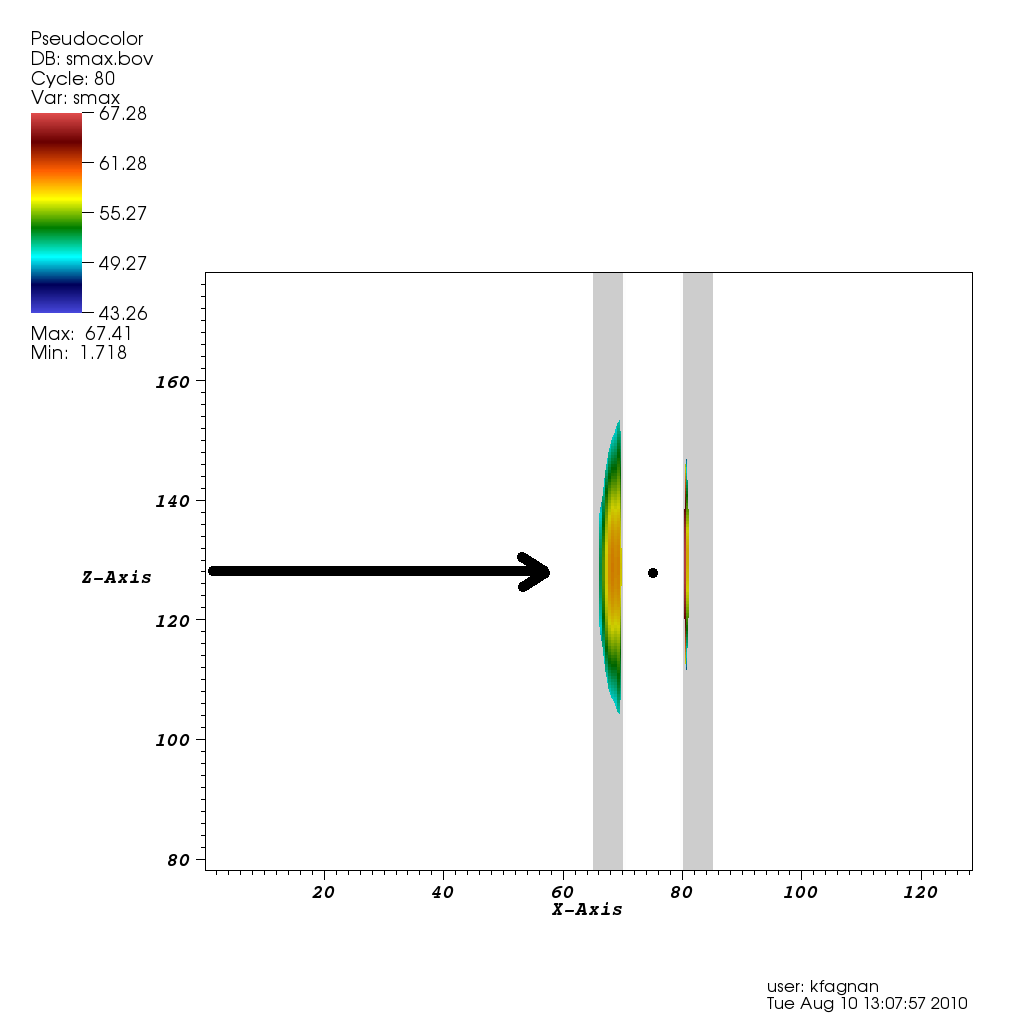}
\caption{An additional interface calculation showing the more realistic treatment of a cylinder with a fluid-filled cavity.  a) is a slice along z=0 showing the concentration of stress in the front of the idealized bone, with additional smaller pockets of maximum stress due to reflection in the back half of the bone.  b) is a slice along y=0 which further demonstrates the stress concentration in the first half of the bone.  The arrow in b indicates the direction of ESWT wave propagation and the dot indicates the location of F2.\label{fig:bone_canal}}
\end{center}

\end{figure}

In Figure \ref{fig:bone_canal} we show two-dimensional slices of a calculation with a more realistic, but still idealized, long bone geometry.  In a the shaft of a long bone, there is a marrow-filled canal running through the center.  Marrow is typically modeled as a viscoelastic material \cite{martin_burr_sharkey}, but for this first approximation we just used a fluid-filled canal.  The impedance difference in the two materials is similar and therefore illustrates the behavior that we are interested in, the change in maximal stress deposition.  In contrast to the solid cylinder above, the contours of maximal stress are concentrated in the front side of the hollow cylinder.  Figure \ref{fig:bone_canal} b) shows that there are also two regions of additional stress concentration in the backside of the hollow cylinder.  This example highlights the importance in understanding where these impedance jumps occur in order to optimally treat the patient.

\subsection{Heterotopic Ossification}
\label{sec:ho}
\begin{figure}[h!]
\begin{center}
\includegraphics[height=2in]{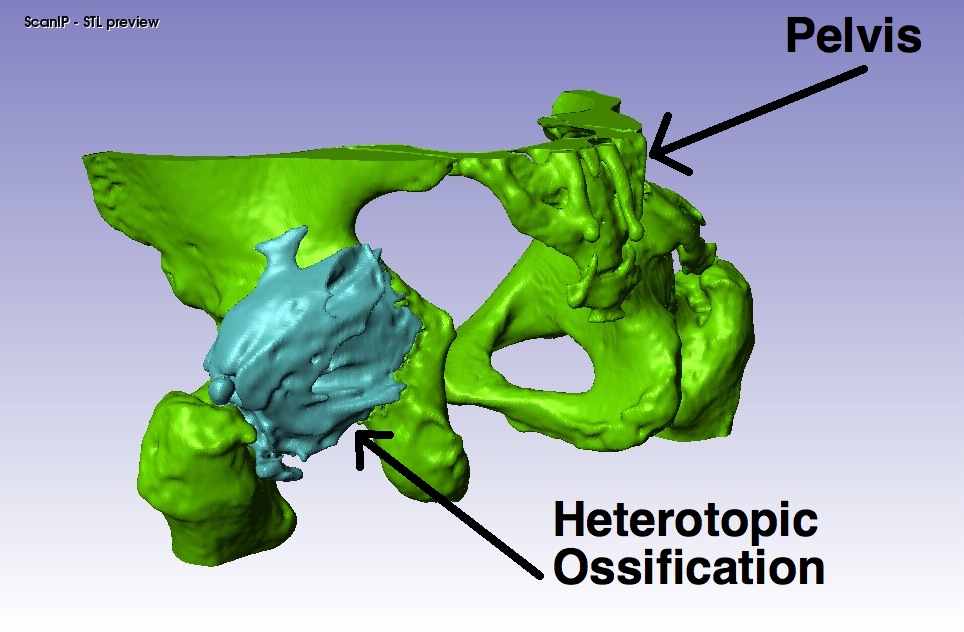}\hspace{5mm}
\includegraphics[height=2in]{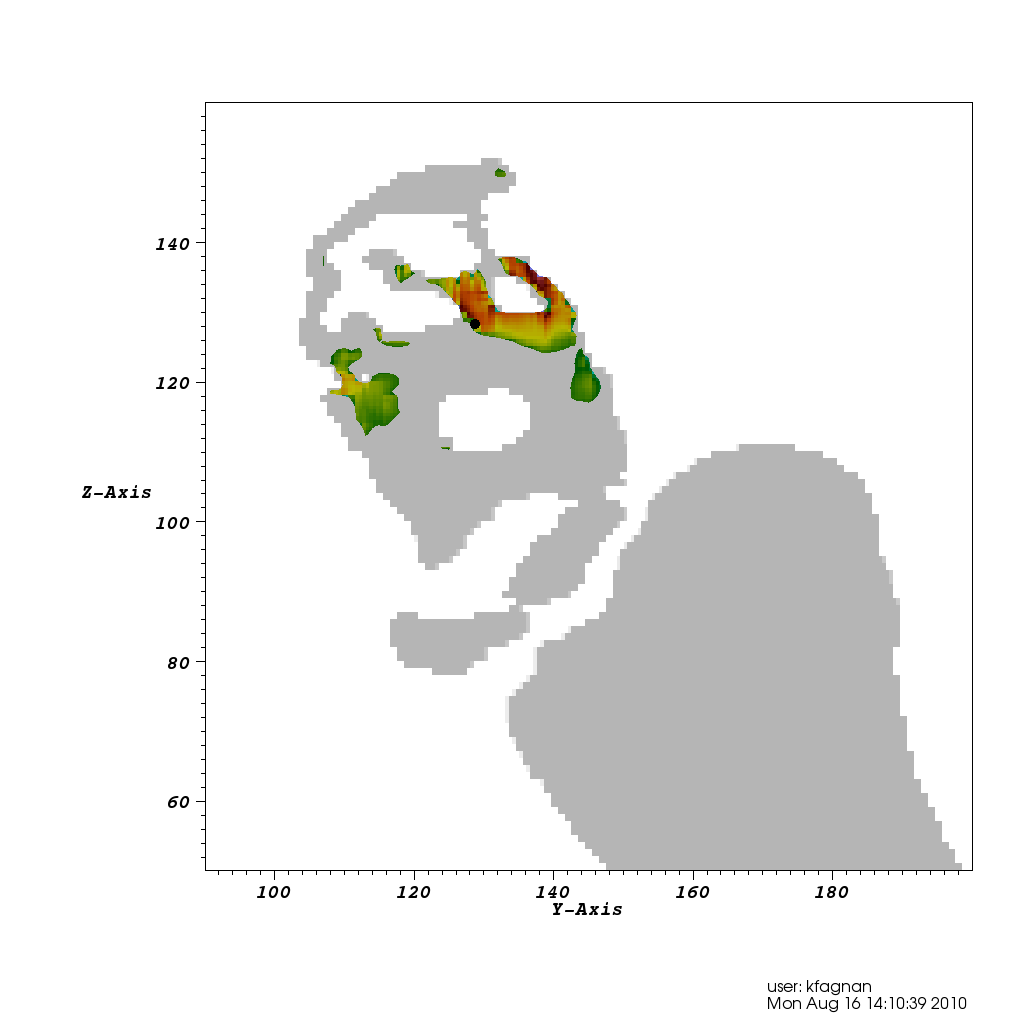}
\caption{a) The 3D CT patient data illustrating the heterotopic ossification (blue) attached to the right hip joint (green).  b) A slice at x=115 of the 2D calculation shows how the pockets of fluid lead to stress concentration in the substructure of the ossification, the dot indicates the location of F2 and the direction of treatment is into the page.}
\label{fig:real_ho}
\end{center}
\end{figure}

A heterotopic ossification (HO) is a growth of bone-like 
material in soft tissue.  HOs often grow spontaneously in tissue that has been traumatized due to injury or amputation.
An example of an HO is shown in Figure \ref{fig:real_ho} a), which shows the pelvis and an HO, using data extracted from a patient's CT scan.  In this case, the HO has grown around the right hip joint and is inhibiting the patient's range of motion.  The goal of the HO treatment is not to pulverize the ossification, but to break up the adhesion between the HO and the joint, in order to restore the patient's range of motion.  There is no clear division between the HO and bone in the CT scan because both are composed of materials that have similar densities.  However, the HO does not have the same woven structure that is present in bone, so the two will likely have different material properties, even though the densities are similar.  This similarity means that we are uncertain as to how strong the connection or adhesion is between the HO and the bone, which will directly impact the number of shocks needed to restore the patient's range of motion.

We are able to use our model to investigate the effect of the angle of treatment on the observed stresses in the region near the HO.  Since the composition and material 
properties of the ossification are not well understood, we can also use the model to
vary the material properties of the ossification and investigate the sensitivity of the results to these parameters.  We found that both the strength of the connection between the HO and bone, as well as the composition of the HO, had a significant effect on the location of maximum 
stress in the object \cite{kfagnan_mchang_ho}.  

It is challenging to infer anything meaningful from the images in the full three-dimensional calculation in Figure \ref{fig:real_ho} a), so we have also included a two-dimensional slice of the maximum shear in Figure \ref{fig:real_ho} b).  Here the gray regions represent the bone-like HO material and we assume any gaps are filled with fluid.  It is clear that the interior of the ossification is complex and contains may fluid-filled pockets that affect, in this case, the location of the maximum shear stress.  

\begin{figure}[h!]
\begin{center}
\includegraphics[height=2.5in]{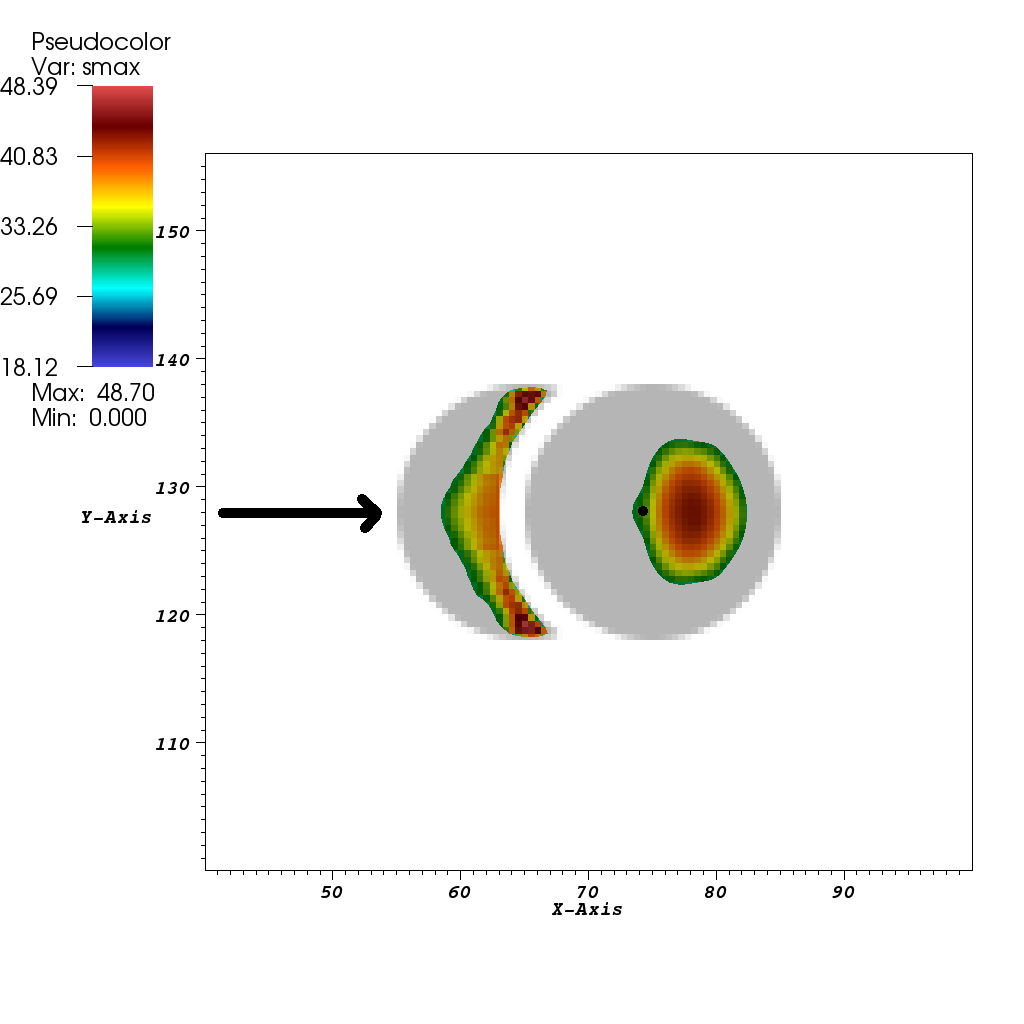}\hspace{5mm}
\includegraphics[height=2.5in]{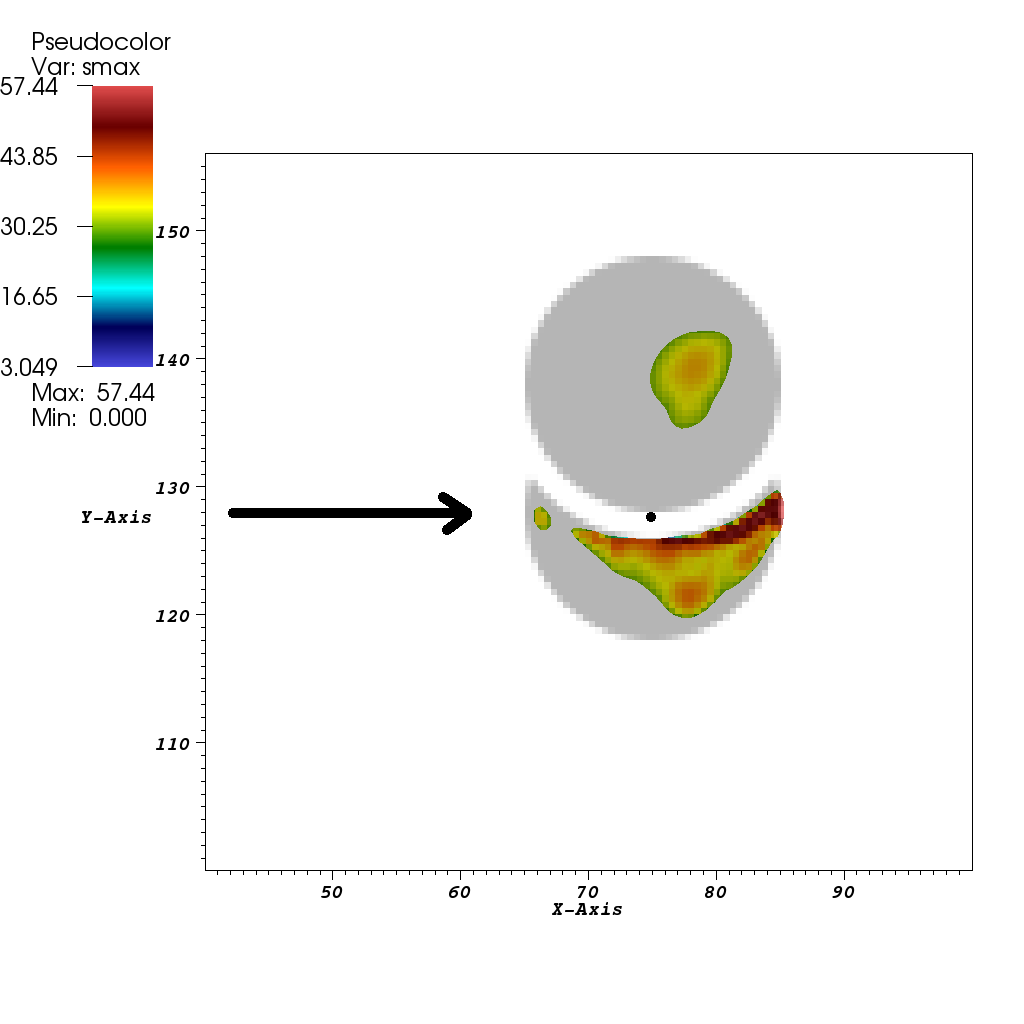}
\caption{Calculations for an idealized ossification that demonstrates the difference in shear stress deposition when treating the HO from different directions.  Since the goal is to disrupt the HO at the interface between the bone and the HO, b) indicates that it would be optimal send shock waves parallel to the break, instead of perpendicular to it.  The arrows indicate the direction of ESWT propagation and the dots indicate the location of F2.}
\label{fig:broken_ho_shear}
\end{center}
\end{figure}

Given the complex nature of the HO and subsequent difficulty interpreting the three-dimensional results, we have also used an idealized ossification to investigate some facets of the treatment.  One example of this is shown in Figure \ref{fig:broken_ho_shear}, where we have simulated a case where the ossification (the crescent in the two-dimensional images) is not strongly attached to the bone (the cylinder) and calculated the maximal shear stress as a result of two different treatment angles.  Figure \ref{fig:broken_ho_shear} a) illustrates the result when the ESWT device is aimed orthogonal to the gap between the HO and the cylinder.   Figure \ref{fig:broken_ho_shear} b) is the result when the device is aimed so the shockwaves propagate parallel to the gap between the HO and bone.  It has been indicated that maximum shear stress is important in causing the HO to break, so it is desirable to deposit the maximum amount of shear as close to the HO/bone interface as possible.  In this case, it is better to treat the HO in the direction indicated in Figure \ref{fig:broken_ho_shear} b), since the shear stress is concentrated along the gap.  

According to our computational results, the pockets of fluid within an HO and strength of adhesion to the bone surface will affect the stress deposition and therefore the location of the eventual break in the ossification.  Further investigation is required to be conclusive, but our results indicate that if the fluid pockets are in the propagation path of the shock wave, they may cause the maximum stresses to occur away from the adhesion site, making the treatment less effective.  In reality, the composition of the HO is unknown and we do not have a good characterization for the material properties of the ossifications.  However, the strong impedance mismatch between fluid and bone, as well as the inability of the fluid to support shear stress, indicate that the presence of fluid-filled pockets will have an effect on the stress deposition.   
We should note here that our modeling work does not take into account the propagation of successive shocks or failure models within the material, which should ultimately be incorporated in order to  the determine the optimal treatment.  This is an area for future work.

\ignore{ 
\begin{figure}[ht]
\begin{center}
a)\includegraphics[scale=.25]{hip_comp_direct.png}
b)\includegraphics[scale=.25]{hip_tens_direct_ydir.png}
\caption{\alert{these will be a 3d object} Isosurface of a) Maximum Compression and b) Maximum Tension in the HO shot directly.}
\label{fig:3d_ho}
\end{center}
\end{figure}
\alert{These will be a 3D embedded object, once I have figured out how to plot the stresses on the same 
figure and have the shading be different - I haven't yet determined how to do that.}
}

\section{Conclusion}
In this paper we have proposed a new model for ESWT.  We have demonstrated that the Tait 
equation of state is sufficient for the pressures that arise in ESWT.  We have shown that the fluid and 
solid can be modeled with the same set of Lagrangian equations since the particle displacements are 
small.  This approach allowed us to utilized existing numerical methodology, consisting of 
high-resolution shock-capturing methods together with Adaptive Mesh Refinement, to efficiently 
calculate solutions to these equations for a variety of idealized biological problems.  We have 
also demonstrated that we can effectively handle interfaces between different materials on Cartesian 
grids.  Using this methodology we were able to explore, even in geometrically complicated structures, 
how the interfaces between the fluid and solid materials  affect the distribution of maximal 
stress in several problems of clinical interest.

Maximizing stress in specific regions seems important in both the
healing and destruction of biological tissues.  Shear stress is
thought to be play a role in the stimulation of biological
tissues \cite{eswt_vegf,frangos,freund,goodship1985,claes1999,lacroix2002,park1998}.  
Mechanical loading is thought 
to play a role in the formation of bone tissue, and as discussed in \Sec{intro}, 
shear and compressive displacements generated by loading influence bone healing 
\cite{park1998,prendergast1997,lacroix2002,claes1999,isaksson2006,carter1998,goodship1985}.  
Shear stress is also important in predicting the break up of kidney stones \cite{bailey_oleg}.  Tension
leads to the formation of cavitation bubbles
that lead to damage and crack formation in kidney stones when
subjected to shock wave lithotripsy.  These cracks are thought to
be one mechanism that initiates kidney stone comminution, or the pulverization of the stone.  The role
of cavitation in the treatment of heterotopic ossifications is not
yet well understood, but given the similarity of HO treatments to
lithotripsy, it is reasonable to hypothesize that cavitation will
impact the treatment.

The model we have developed has been used to investigate idealized non-unions and 
heterotopic ossifications, and we have shown a few examples to illustrate this.  
A broader range of calculations are available in \cite{fagnan:phd}.  These specific examples are
 being studied in more detail and papers with additional work on these problems are now in 
 preparation \cite{kfagnan_mchang_ho, amath_apl_nonunion}.
  
The focus of this paper has been the effect that material interfaces between
tissue and bone have on the transmission, reflection, and 
focusing of the shock wave.  Very simple models have been used for the
material on each side of the interface: compressible fluid with a Tait
equation of state in the tissue and linear isotropic elasticity in the bone.
We believe that this level of macroscopic modeling can already reveal
interesting features of the stress that may be clinically important.  In
particular, focusing may occur in regions displaced from where it would be
observed in pure water, and mode conversion at an interface can generate
shear waves in the bone that are not present in the focusing shock
wave in fluid.

To consider the effect of stress on individual osteocytes,
a much more detailed model would be necessary that is beyond the scope
of this work.  In particular, 
this would require modeling the microscale fluid-filled canaliculi
within the bone through which the osteocyte processes extend. Work is
currently underway in this direction, and also on intermediate levels of
modeling in which the bone is modeled as an orthotropic poroelastic
material.  These equations can be solved with essentially the same high
resolution finite volume methods used here, after implementing a more
complicated Riemann solver \cite{lemoine:pc}, and with the same software for
adaptive mesh refinement.
Another possible extension is to
investigate viscoelastic tissue models that may be superior to the
Tait equation for water that is currently used.

\section*{Acknowledgements}
This work was supported in part by NIH Grant 5R01AR53652-2,
NSF Grants DMS-0609661 and DMS-0914942, and the Founders Term Professorship
in Applied Mathematics at the University of Washington.
The authors would like to thank Donna Calhoun and the ANAG 
group at Lawrence Berkeley Laboratory for their assistance with the ChomboClaw calculations.
This research was supported in part by the National Science Foundation through 
TeraGrid resources provided by TACC under grant number TG-DMS090036T.

\bibliographystyle{plain}
\bibliography{references}

\begin{thebibliography}{10}

\bibitem{augat2001}
P~Augat, J~Merk, S~Wolf, and LE~Claes.
\newblock Mechanical stimulation by external application of cyclic tensile
  strains does not effectively enhance bone healing.
\newblock {\em Journal of Orthopaedic Trauma}, 15:54--60, 2001.

\bibitem{cleveland_averkiou}
M.A. Averkiou and R.O. Cleveland.
\newblock Modeling of an electrohydraulic lithotripter with the {KZK} equation.
\newblock {\em Journal for the Acoustical Society of America}, 106(1):102--112,
  1999.

\bibitem{eswt_evidence_2010}
MD~B.~Zelle and et. al.
\newblock Extracorporeal shock wave therapy: Current evidence.
\newblock {\em Journal of Orthopaedic Trauma}, 24(3):s66--s70, 2010.

\bibitem{db-rjl-sm-jr:vcflux}
D.~Bale, R.~J. LeVeque, S.~Mitran, and J.~A. Rossmanith.
\newblock A wave-propagation method for conservation laws and balance laws with
  spatially varying flux functions.
\newblock {\em SIAM J. Sci. Comput.}, 24:955--978, 2002.

\bibitem{mjb-ol:amr}
M.~Berger and J.~Oliger.
\newblock Adaptive mesh refinement for hyperbolic partial differential
  equations.
\newblock {\em J. Comput. Phys.}, 53:484--512, 1984.

\bibitem{mjb-col:amr}
M.~J. Berger and P.~Colella.
\newblock Local adaptive mesh refinement for shock hydrodynamics.
\newblock {\em J. Comput. Phys.}, 82:64--84, 1989.

\bibitem{mjb-rjl:amrclaw}
M.~J. Berger and R.~J. LeVeque.
\newblock Adaptive mesh refinement using wave-propagation algorithms for
  hyperbolic systems.
\newblock {\em SIAM J. Numer. Anal.}, 35:2298--2316, 1998.

\bibitem{mjb-rig:cluster}
M.~J. Berger and I.~Rigoutsos.
\newblock An algorithm for point clustering and grid generation.
\newblock {\em IEEE Trans. Sys. Man \& Cyber.}, 21:1278--1286, 1991.

\bibitem{eswt_bone}
R.~Biedermann, A.~Martin, G.~Handle, T.~Auckenthaler, C.~Bach, and M.~Krismer.
\newblock Extracorporeal shock waves in the treatment of nonunions.
\newblock {\em Journal of Trauma}, 54(5):936--42, 2003.

\bibitem{chomboclaw}
D.~A. Calhoun, P.~Colella, and R.~J. LeVeque.
\newblock {\sc chombo-claw} software.
\newblock \\ {\tt http://www.amath.washington.edu/~calhoun/demos/ChomboClaw}.

\bibitem{carter1998}
D.R. Carter, G.S. Beaupre, N.J. Giori, and J.A. Helms.
\newblock Mechanobiology of skeletal regeneration.
\newblock {\em Clinical Orthopaedics and Related Research}, 355(Suppl):S41--55,
  1998.

\bibitem{christopher_hm3}
T.~Christopher.
\newblock Modeling the {Dornier HM3} lithotripter.
\newblock {\em Journal of the Acoustical Society of America}, 96(5):3088--3095,
  1994.

\bibitem{claes1999}
L.E. Claes and C.A. Heigele.
\newblock Magnitudes of local stress and strain along osseous surfaces predict
  the course and type of fracture-healing.
\newblock {\em The Journal of Biomechanics}, 32:255--266, 1999.

\bibitem{claes1995}
L.E. Claes, H.J. Wilke, P.~Augat, S.~Rubenacker, and K.J. Margevicius.
\newblock Effect of dynamization on gap healing of diaphyseal fractures under
  external fixation.
\newblock {\em Cinical Biomechanics}, 10:227--234, 1995.

\bibitem{oleg_cleveland}
R.O. Cleveland and O.~Sapozhnikov.
\newblock Modeling elastic wave propagation in kidney stones with application
  to shock wave lithotripsy.
\newblock {\em Journal of the Acoustical Society of America},
  118(4):2667--2676, 2005.

\bibitem{colella_communication}
P.~Colella.
\newblock Personal communication, 2009.

\bibitem{coleman}
A.J. Coleman, J.~Saunders, R.~Preston, and D.~Bacon.
\newblock Pressure waveforms generated by a {Dornier} extra-corporeal
  shock-wave lithotripter.
\newblock {\em Ultrasound in Medicine and Biology}, 13:651--657, 1987.

\bibitem{lacroix2002}
P.J.~Prendergast D.~Lacroix.
\newblock A mechano-regulation model for tissue differentiation during
  fracture-healing: analysis of gap size and loading.
\newblock {\em The Journal of Biomechanics}, 35:1163--1171, 2002.

\bibitem{clawpack}
R.J.~LeVeque et. al.
\newblock \url{http://www.clawpack.org}, 2006.

\bibitem{fagnan:phd}
K.~Fagnan.
\newblock {\em High-resolution finite volume methods for extracorporeal shock
  wave therapy}.
\newblock PhD thesis, University of Washington, 2010.

\bibitem{kfagnan_mchang_ho}
K.~Fagnan, M.~Chang, and R.J. LeVeque.
\newblock Computational investigation of {ESWT} for treatment of heterotopic
  ossifications.
\newblock {\em in preparation}, 2010.

\bibitem{kfagnan_hyp06}
K.M. Fagnan, R.J. LeVeque, T.J. Matula, and B.~MacConaghy.
\newblock High-resolution finite volume methods for extracorporeal shock wave
  therapy.
\newblock In S{ylvie} Benzoni-Gavage and D{enis} Serre, editors, {\em
  Hyperbolic Problems: Theory, Numerics, Applications}, pages 503--510.
  Springer, 2006.

\bibitem{freund}
J.~Freund, T.~Colonius, and A.~Evan.
\newblock A cumulative shear mechanism for tissue damage initiation in
  shock-wave lithotripsy.
\newblock {\em Ultrasound in Medicine and Biology}, 33:1495--1503, 2007.

\bibitem{fung}
Y.C. Fung.
\newblock {\em Biomechanics: Mechanical Properties of Living Tissues}.
\newblock Springer, 1993.

\bibitem{visc_bone}
E.~Garner, R.~Lakes, T.~Lee, C.~Swan, and R.~Brand.
\newblock Viscoelastic dissipation in compact bone: Implications for
  stress-induced fluid flow in bone.
\newblock {\em Journal of Biomechanical Engineering}, 122:166--73, 2000.

\bibitem{goodship1985}
A.E. Goodship and J.~Kenwright.
\newblock The influence of induced micromovement on the healing of experimental
  tibial fractures.
\newblock {\em Journal of Bone and Joint Surgery British Volume}, 67:650--655,
  1985.

\bibitem{hamilton}
M.~Hamilton.
\newblock Transient axial solution for the reflection of a spherical wave from
  a concave ellipsoidal mirror.
\newblock {\em Journal of the Acoustical Society of America}, 93(3):1256--1266,
  1993.

\bibitem{frangos}
M.V. Hillsley and J.A. Frangos.
\newblock Review: Bone tissue engineering: The role of interstitial fluid flow.
\newblock {\em Biotechnology and Bioengineering}, 43:573--581, 1994.

\bibitem{huang2010}
C.~Huang and R.~Ogawa.
\newblock Mechanotransduction in bone repair and regeneration.
\newblock {\em The FASEB Journal}, 23:3625--3632, 2010.

\bibitem{isaksson2006}
H~Isaksson, W~Wilson, CC~van Donkelaar, R~Huiskes, and K~Ito.
\newblock Comparison of biophysical stimuli for mechano-regulation of tissue
  differentiation during fracture-healing.
\newblock {\em The Journal of Biomechanics}, 39:1507--1516, 2006.

\bibitem{ivings_toro}
M.J. Ivings, D.M. Causon, and E.F. Toro.
\newblock On {Riemann} solvers for compressible liquids.
\newblock {\em International Journal for Numerical Methods in Fluids},
  28:395--418, 1998.

\bibitem{orthotripsy_ogden}
MD~J.A.~Ogden, MD~R.~G.~Alvarez, MD~R.~Levitt, and RN~M.~Marlow.
\newblock Shock wave therapy (orthotripsy) in musculoskeletal disorders.
\newblock {\em Clinical Orthopaedics and Related Research}, 387:22--40, 2001.

\bibitem{keaveny}
T.~Keaveny, X.E. Guo, E.F. Wachtel, T.A. McMahon, and W.C. Hayes.
\newblock Trabecular bone exhibits fully linear elastic behavior and yields at
  low strains.
\newblock {\em Journal of Biomechanics}, 27, 1994.

\bibitem{amath_apl_nonunion}
K.Fagnan, T.~Matula, and R.J. LeVeque.
\newblock Optimized {ESWT} treatment for nonunions.
\newblock {\em in preparation}, 2010.

\bibitem{chombo}
Applied Numerical Algorithms Group (ANAG) Lawrence Berkeley~National
  Laboratory.
\newblock \url{https://seesar.lbl.gov/ANAG/chombo/}, 2009.

\bibitem{jol-rjl:3d}
J.~O. Langseth and R.~J. LeVeque.
\newblock A wave-propagation method for three-dimensional hyperbolic
  conservation laws.
\newblock {\em J. Comput. Phys.}, 165:126--166, 2000.

\bibitem{lemoine:pc}
G.~Lemoine.
\newblock in preparation.

\bibitem{rjl:wpalg}
R.~J. LeVeque.
\newblock Wave propagation algorithms for multi-dimensional hyperbolic systems.
\newblock {\em J. Comput. Phys.}, 131:327--353, 1997.

\bibitem{claw.org.url}
R.~J. LeVeque, M.~J. Berger, et~al.
\newblock {\sc clawpack} software.
\newblock \href{http://www.clawpack.org}{www.clawpack.org}.

\bibitem{rjl_book}
R{andall}.~J. LeVeque.
\newblock {\em Finite Volume Methods for Hyperbolic Problems}.
\newblock Cambridge University Press, 2002.

\bibitem{rjl_nonlinear}
R.J. LeVeque.
\newblock Finite-volume methods for non-linear elasticity in heterogeneous
  media.
\newblock {\em International Journal for Numerical Methods in Fluids},
  40(1-2):93--104, 2001.

\bibitem{martin_burr_sharkey}
R.~Bruce Martin, David~B. Burr, and Neil~A. Sharkey.
\newblock {\em Skeletal Tissue Mechanics}.
\newblock Springer, 1998.

\bibitem{matula_direct}
T.J. Matula, P.R. Hilmo, and M.R. Bailey.
\newblock A suppressor to prevent direct wave-induced cavitation in shock wave
  therapy devices.
\newblock {\em Journal of the Acoustical Society of America}, 118(1):178--185,
  2005.

\bibitem{morgan2008}
E.~F. Morgan, R.~E. Gleason, L.~N.~M. Hayward, P.L. Leong, and K.T.
  Salisbury-Paolomares.
\newblock Mechanotransduction and fracture repair.
\newblock {\em Journal of Bone and Joint Surgery American Volume}, 90(Suppl
  1):25--30, 2008.

\bibitem{eswt_shoulder}
G.~Mouzopoulos, M.~Stamatakos, D.~Mouzopoulos, and M.~Tzurbakis.
\newblock Extracorporeal shock wave treatment for shoulder calcific tendonitis:
  a systematic review.
\newblock {\em Skeletal Radiology}, 36(9):803--811, 2008.

\bibitem{nakahara}
M.~Nakahara, K.~Nagayama, and Y.~Mori.
\newblock Shockwave dynamics of high pressure pulse in water and other
  biological materials based on {Hugoniot} data.
\newblock {\em Japanese Journal of Applied Physics}, 47, 2008.

\bibitem{ogden}
J.A. Ogden, A.~Toth-Kischkat, and R.~Schultheiss.
\newblock Principles of shock wave therapy.
\newblock {\em Clinical Orthopaedics and Related Research}, 387:8--17, 2001.

\bibitem{park1998}
SH~Park, K~O'Connor, H~McKellop, and A~Sarmiento.
\newblock The influence of active shear or compressive motion on
  fracture-healing.
\newblock {\em The Journal of Bone and Joint Surgery American Volume},
  80:868--878, 1998.

\bibitem{prendergast1997}
P.J. Prendergast, R.~Huiskes, and K.~Soballe.
\newblock Biophysical stimuli on cells during tissue differentiation at implant
  interfaces.
\newblock {\em The Journal of Biomechanics}, 30:539--548, 1997.

\bibitem{robling2002}
AG~Robling, FM~Hinant, DB~Burr, and CH~Turner.
\newblock Improved bone structure and strength after long-term mechanical
  loading is greatest if loading is separated into short bouts.
\newblock {\em Journal of Bone Mineral Research}, 17:1545--1554, 2002.

\bibitem{roblingcyclic2002}
AG~Robling, FM~Hinant, DB~Burr, and CH~Turner.
\newblock Shorter, more frequent mechanical loading sessions enhance bone mass.
\newblock {\em Medicine and Science in Sports and Exercise}, 34:196--202, 2002.

\bibitem{saito}
T.~Saito, M.~Marumoto, H.~Yamashita, S.H.R. Hosseini, A.~Nakagawa, T.~Hirano,
  and K.~Takayama.
\newblock Experimental and numerical studies of underwater shock wave
  attenuation.
\newblock {\em Shock Waves}, 13:139--148, 2003.

\bibitem{oleg_bailey}
O.~Sapozhnikov, M.~Bailey, and R.O. Cleveland.
\newblock The role of shear and longitudinal waves in the kidney stone
  comminution by a lithotripter shock pulse (a).
\newblock {\em Journal of the Acoustical Society of America}, 115:2562--2562,
  2004.

\bibitem{bailey_oleg}
O.~Sapozhnikov, A.D. Maxwell, B.~MacConaghy, and M.~Bailey.
\newblock A mechanistic analysis of stone fracture in lithotripsy.
\newblock {\em Journal of the Acoustical Society of America},
  121(2):1190--1202, 2007.

\bibitem{saxon2005}
LK~Saxon, AG~Robling, I~Alam, and CH~Turner.
\newblock Mechanosensitivity of the rat skeleton decreases after a long period
  of loading, but is improved with time off.
\newblock {\em Bone}, 36:454--464, 2005.

\bibitem{tanguay}
M.~Tanguay.
\newblock {\em Computation of bubbly cavitating flow in shock wave
  lithotripsy}.
\newblock PhD thesis, California Institute of Technology, 2004.

\bibitem{orthotropic}
W.R. Taylor, E.~Roland, H.~Ploeg, D.~Hertig, R.~Klabunde, M.D. Warner, M.C.
  Hobatho, L.~Rakotomanana, and S.E. Clift.
\newblock Determination of orthotropic bone elastic constants using {FEA} and
  modal analysis.
\newblock {\em Journal of Biomechanics}, 35:767--73, 2002.

\bibitem{turner1998}
CH~Turner and FM~Pavalko.
\newblock Mechanotransduction and functional response of the skeleton to
  physical stress: the mechanisms and mechanics of bone adaptation.
\newblock {\em Journal of Orthopaedic Science}, 3:346--355, 1998.

\bibitem{valchanou}
VD~Valchanou and P~Michailov.
\newblock High energy shock waves in the treatment of delayed and nonunion of
  fractures.
\newblock {\em International Orthopaedics}, 15(3):181--184, 1991.

\bibitem{eswt_vegf}
C.-J. Wang, K.-E. Huang, Y.-C. Sun, Y.-J. Yang, J.-Y. Ko, L.-H. Weng, and F.-S.
  Wang.
\newblock Vegf modulates angiogenesis and osteogenesis in shockwave-promoted
  fracture healing in rabbits.
\newblock {\em Journal of Surgical Research}, In Press, Corrected Proof:--,
  2010.

\bibitem{cjwang_hip}
C.-J. Wang, F.-S. Wang, J.-Y. Ko, H.-Y. Huang, C.-J. Chen, Y.-C. Sun, and Y.-J.
  Yang.
\newblock Extracorporeal shockwave therapy shows regeneration in hip necrosis.
\newblock {\em Rheumatology}, 47(4):542--546, 2008.

\bibitem{weinbaum1994}
S~Weinbaum, SC~Cowin, and Y~Zeng.
\newblock A model for the excitation of osteocytes by mechanical
  loading-induced bone fluid shear stress.
\newblock {\em Journal of Biomechanics}, 27:339--360, 1994.

\bibitem{wolf2001}
S~Wolf, P~Augat, K~Eckert-Hubner, A~Laule, GD~Krischak, and LE~Claes.
\newblock Effects of high-frequency, low-magnitude mechanical stimulus on bone
  healing.
\newblock {\em Clinical Orthopaedics and Related Research}, 385:192--198, 2001.

\bibitem{zelle2010}
BA~Zelle, H~Gollwitzer, M~Zlowodzki, and V~Burhen.
\newblock Extracorporeal shock wave therapy: Current evidence.
\newblock {\em Journal of Orthopaedic Trauma}, 24(3 (Suppl)):S66--S70, 2010.

\end{thebibliography}

\end{document}